\documentclass[aps,prd,reprint,twocolumn,groupedaddress]{revtex4-1}

\usepackage{algorithm}
\usepackage{algpseudocode}

\usepackage{graphicx} % Include figure files
\usepackage{dcolumn}  % Align table columns on decimal point
\usepackage{bm}       % Bold math
\usepackage{amsmath}  % Enhanced math support
\usepackage{amssymb}  % Additional math symbols
\usepackage{natbib}
\usepackage{dsfont}
\usepackage{bbold}
\usepackage{nicematrix}
\usepackage[colorlinks=true, urlcolor=violet, linkcolor=blue, citecolor=red, hyperindex=true, linktocpage=true, draft=false]{hyperref} 

\usepackage{tikz-cd}
\usepackage{xcolor}
\usepackage{multirow}
\usepackage{tablefootnote}
\usepackage{natbib}
\usetikzlibrary{positioning}

\definecolor{Pank}{rgb}{1.0,0.1,0.5}

\usepackage{cancel}

\usepackage{thm-restate}

\newtheorem{lemX}{Lemma}

\definecolor{darkblue}{rgb}{0,0,.6}
\definecolor{orange}{rgb}{1,0.5,0}

\usepackage{comment}

\usepackage{wasysym,mdframed,tikz,leftindex,booktabs,caption}
\usetikzlibrary{patterns,arrows}
\NiceMatrixOptions{cell-space-limits = 2pt}

\newcommand{\I}{\mathrm{I}} %I ~ Identity operator
\newcommand{\J}{\mathrm{J}}
\newcommand{\T}{\mathrm{T}}
\newcommand{\G}{\mathrm{G}}
\newcommand{\F}{\mathrm{F}}
\newcommand{\M}{\mathrm{M}}

\definecolor{nicegreen}{rgb}{0.13, 0.55, 0.13}

\begin{document}

% Title of paper
\title{Fast and Flexible Quantum-Inspired PDE Solvers with Data Integration}

% Authors
\author{Lucas Arenstein$^{1}$}\email{lsa@di.ku.dk}
\author{Martin Mikkelsen$^1$}
\author{Michael Kastoryano$^{1,2}$}

\affiliation{$^1$Department of Computer Science, University of Copenhagen, Denmark}
\affiliation{$^2$AWS Center for Quantum Computing, Pasadena, CA, USA}

\date{\today}

\begin{abstract}
Accurately solving high-dimensional partial differential equations (PDEs) remains a central challenge in computational mathematics. Traditional numerical methods, while effective in low-dimensional settings or on coarse grids, often struggle to deliver the precision required in practical applications. Recent machine learning-based approaches offer flexibility but frequently fall short in terms of accuracy and reliability, particularly in industrial contexts. In this work, we explore a quantum-inspired method based on quantized tensor trains (QTT), enabling efficient and accurate solutions to PDEs in a variety of challenging scenarios. Through several representative examples, we demonstrate that the QTT approach can achieve logarithmic scaling in both memory and computational cost for linear and nonlinear PDEs. Additionally, we introduce a novel technique for data-driven learning within the quantum-inspired framework, combining the adaptability of neural networks with enhanced accuracy and reduced training time.
\end{abstract}


\maketitle

% Main text
\section{Introduction}
One of the fundamental challenges in computational mathematics is the efficient solution of partial differential equations (PDEs), especially as the dimensionality of the problem increases or when very fine grids are required.  High-dimensional PDEs are prevalent in various fields, including quantum mechanics, finance, and fluid dynamics, making it imperative to develop methods that can tackle this complexity without prohibitive computational costs.

Finite difference and finite element methods remain widely used due to their generality and ease of implementation, particularly when combined with fast iterative solvers such as (algebraic) multigrid methods, which exploit hierarchical coarse-to-fine structures to accelerate convergence \cite{brandt1977multi,trottenberg2001multigrid}. In parallel, spectral methods have demonstrated exceptional performance for smooth problems, offering exponential convergence rates by projecting solutions onto global basis functions, such as Fourier or Chebyshev polynomials \cite{gottlieb1977numerical,boyd2001chebyshev}. These techniques represent the state of the art in many classical applications, yet they often face limitations in problems involving complex geometry, sharp interfaces, or multiscale phenomena, where maintaining both accuracy and computational efficiency becomes increasingly challenging.

In recent years, machine learning has emerged as a promising framework for solving differential equations, offering new paradigms that complement or bypass traditional discretization techniques. Physics-informed neural networks (PINNs) encode the differential equation and boundary conditions directly into the loss function of a neural network, enabling mesh-free approximation of solutions in complex geometries and high-dimensional settings \cite{raissi2019physics,berg2019data}. PINNs have since been extended to handle stiff systems, inverse problems \cite{long2019pde}, and fractional differential equations \cite{pang2020npinns}, highlighting their flexibility. In parallel, a family of approaches known as neural operators has aimed to learn mappings from function spaces to function spaces—i.e., the solution operator itself—enabling rapid inference across parametric PDE families. These include DeepONets \cite{lu2021learning}, Fourier Neural Operators (FNOs) \cite{li2020multipole}, and the various recent Galerkin, Green’s function, and transformer-based operator networks. These methods have shown strong empirical performance in surrogate modeling and uncertainty quantification, and they offer certain advantages in terms of generalization and adaptability. However, despite their promise, current neural PDE solvers often lack the precision and reliability of classical numerical methods, particularly for stiff, chaotic, or highly oscillatory systems, and they typically require extensive computational resources and time for training. Bridging this gap remains an active area of research, combining insights from numerical analysis, deep learning, and applied mathematics.

In this work, we aim to circumvent the curse of dimensionality by adopting a quantum-inspired approach rooted in methods developed within the quantum chemistry and condensed matter physics communities \cite{White}. Specifically, we employ the quantized tensor train (QTT) formalism—an efficient representation of multi-dimensional functions that leverages tensor decomposition techniques from quantum information theory. The QTT format enables compression and manipulation of large-scale data by exploiting low-rank structure, thereby significantly reducing computational and memory demands when solving PDEs. In this context, a common challenge arises from the trade-off between computational efficiency and solution accuracy: achieving high accuracy often requires fine discretization grids to resolve small-scale or highly oscillatory features, but such discretizations lead to prohibitively large linear systems that strain both time and memory resources.

QTT-based algorithms provide a robust framework for addressing the computational challenges posed by high-dimensional and multiscale problems. By encoding each spatial variable using binary quantization and expressing the overall function as a tensor train (TT), QTT representations efficiently capture small-scale structure and oscillatory behavior. This hierarchical format significantly reduces the number of degrees of freedom, enabling scalable computations in extremely high-dimensional settings. The QTT methodology has been developed through foundational contributions from both the quantum physics and numerical mathematics communities \cite{Oseledets2009, Oseledets2010, Oseledets2010a, Oseledets2011}, and has demonstrated notable success in recent applications to multiscale and oscillatory PDEs \cite{Parsimonious, lindsey2023multiscale, Ye_Loureiro_2024, gourianov2022a, gourianov2022quantum, gourianov2025tensor}. In particular, QTT-based solvers have been employed for elliptic problems \cite{qtt_khoromskij}, high-dimensional parabolic equations \cite{richter2021solving}, and integral equations \cite{corona2017tensor}, achieving impressive reductions in computational cost without sacrificing accuracy. Moreover, the QTT format facilitates fast arithmetic with structured matrices—including banded, Toeplitz, and stiffness matrices—thus enabling the construction of efficient iterative solvers, preconditioners, and matrix exponentials. Together, these developments underscore the versatility and power of QTT-based methods as a scalable and accurate approach for the numerical simulation of complex high-dimensional systems.

In this paper, we want to generalize and extend these methods to solve more complex PDEs. Specifically, we address the incorporation of complicated boundary conditions and source terms, which are essential for accurately modeling physical systems but introduce additional layers of complexity in the computational process. By developing novel techniques to include these factors within the QTT framework, we enhance the applicability and robustness of this approach.

Moreover, we provide substantial evidence that both linear and nonlinear PDEs can be effectively solved directly in space-time, without resorting to traditional time-stepping schemes. By treating time as an additional spatial-like dimension, the space-time formulation enables global solution strategies that capture the full temporal evolution of the system in a single iteration. This approach not only facilitates parallelism across time but also avoids the need for sequential progression, which is inherent to time-stepping methods. In particular, the QTT space-time approach offers a highly compressed representation of the full space-time solution, allowing for efficient computations even in high-dimensional settings. Notably, this formulation circumvents the Courant–Friedrichs–Lewy (CFL) condition, which typically constrains the time step size in explicit schemes to ensure numerical stability. By avoiding such stability restrictions and exploiting the logarithmic structure of the QTT format, our method achieves both high accuracy and scalability—reaching, for the first time, an overall complexity of $\mathcal{O}(\log(N T))$ in the number of spatial $N$ and temporal $T$ degrees of freedom.

Finally, we extend the capabilities of QTT methods to incorporate data-driven boundary and initial conditions, a direction widely explored in the PINN community. Our approach begins by integrating data into the solution process through spline fitting and a new QTT interpolation technique \cite{lindsey2023multiscale}, which provide smooth and accurate representations of the input. We then incorporate the data-learned boundary conditions into the QTT pipeline to leverage the computation and memory advantages of our approach. This hybrid strategy leverages the precision and robustness of classical methods while benefiting from the flexibility and data adaptability of modern machine learning approaches. We thus argue that the QTT approach gets the best of both worlds, and should be considered one of the leading numerical tools for solving high dimensional PDEs to high precision. 

\subsection{Outline}
The remainder of this paper is structured as follows. In Section~\ref{sec:TNs}, we introduce key concepts from tensor networks that will be used throughout this work, including: graphical notation, the QTT representation, methods to decompose functions in QTT format, a finite difference scheme in QTT format, and the linear solver ALS. Section~\ref{sec:fd_poisson} presents our QTT-based solver for the Poisson equation, along with benchmark results for 2D and 3D examples, comparing to an algebraic multigrid solver. Section~\ref{sec:fd_burgers} presents our time-stepping and space-time QTT solvers for Burgers’ equation, concluding with a runtime vs. MSE comparison of the two approaches on a specific problem.
Finally, in Section~\ref{sec:data}, we demonstrate how to extend our solver to include learning from data.

\section{Tensor Networks}\label{sec:TNs}
\subsection{The tensor train}
We provide a brief introduction to  tensor network concepts, decompositions, and the notation that will be used throughout this work.
An tensor $\T \in \mathbb{R}^{n_1 \times \cdots \times n_d}$ is a $d$-dimensional array where the indices range from $1$ to $n_k$ in the $k$th \textit{mode}.  An element of the tensor $\T$ is denoted by $\T_{x_1,\hdots,x_d}$ with $1 \leq x_i \leq n_i$ for $i=1,\hdots,d$. Large tensors can be constructed as a network of smaller tensors where common indices are summed over. In this work, we will primarily work with tensor trains, which is a specific tensor network on a line \cite{MPS_Paper,MPSold,tt_ose}. Specifically, the tensor with indices $x_1,...,x_d$ is a tensor train if it can be written as a product of matrices $\{A_j^{x_j}\}_{j=1}^d$ as 
\begin{equation}
T_{x_1,x_2,...,x_d}=A_1^{x_1} A_2^{x_2} \cdots A_d^{x_d}.
\end{equation}
The left and rightmost matrices $A_1^{x_1}$ and $A_d^{x_d}$ are $1\times R_1$ and $R_{d-1}\times 1$ respectively, while all other sets of matrices are of dimensions $R_j\times R_{j+1}$. The $R_j$ are commonly called \textit{bond dimensions} and capture the correlation between the tensor indices along the chain. Note that $A_j$ are 3-tensors of dimension $(R_{j-1},n_j,R_j)$.   

An essential extension of the TT is the matrix product operator (MPO), which describes a factorization of a tensor with $d$ input and $d$ output indices into a line of smaller  (rank-$4$) tensors of dimension $(R_{j-1},n_j,n_j,R_j)$. For TT/MPO,  the uncontracted indices $\{n_j\}$ are called \textit{physical} degrees of freedom, while the contracted indices $\{R_j\}$ are called \textit{virtual} degrees of freedom. Tensor trains serve as compressed representations of large vectors in high-dimensional spaces, while MPOs represent matrices. The essential feature of TTs and MPO, is that certain linear operations, such as matrix-vector multiplication can be performed at the level of individual tensors, hence dramatically reducing the memory and computational footprint.

 \subsection{Graphical notation}
 
 Tensor Network Notation (TNN) provides a powerful tool for visualizing the interactions between tensors in a tensor network. In these visualizations (diagrams), each tensor is represented as a node, with the number of legs corresponding to its dimensions. 
For example:
\begin{itemize}
    \item A matrix \( W \in \mathbb{R}^{m \times n} \) is as a node with two legs: \tikz[baseline=-0.5ex]{
    \node[draw, circle, inner sep=1pt] (tensor) {\scriptsize W};
    \draw (tensor) -- ++(-0.4,0);
    \draw (tensor) -- ++(0.4,0);
}.
    \item A vector \( x \in \mathbb{R}^n \) is a node with a single leg: \tikz[baseline=-0.5ex]{
    %\node[draw, circle, inner sep=0.5pt] (tensor1) {$T$};
    \node[draw, circle, inner sep=0.5pt] (tensor2) at (0.3,0) {x}; % Position tensor2 manually
    %\draw (tensor1) -- (tensor2);
    \draw (tensor2) --(-0.1,0);
}.
\item Matrix-vector multiplication $W x$ is depicted by  \textit{contracting} (summing over) legs of the connected tensors: \tikz[baseline=-0.5ex]{
    \node[draw, circle, inner sep=0.5pt] (tensor1) {\scriptsize W};
    \node[draw, circle, inner sep=0.5pt] (tensor2) at (0.6,0) {x}; % Position tensor2 manually
    \draw (tensor1) -- (tensor2);
    \draw (tensor1) -- ++(-0.3,0);
}
\item Larger tensors $T\in \mathbb{R}^{m\times n\times r}$ have more legs: \tikz[baseline=-0.5ex]{
    \node[draw, circle, inner sep=1pt] (tensor) {\scriptsize T};
    \draw (tensor) -- ++(-0.3,0); % Left leg
    \draw (tensor) -- ++(0.3,0);  % Right leg
    \draw (tensor) -- ++(0,0.3);  % Upper leg
}. 
\end{itemize}
Individual tensor components can be assembled into a tensor network, which can be visualized as a graph: each node represents a tensor, and the connecting lines (or ``legs'') correspond to indices. An edge between two nodes denotes a contracted index (i.e., one over which summation is performed), while legs that are not connected to other nodes represent physical degrees of freedom. In this graphical language, a tensor train is depicted as: 
 \tikz[baseline=-0.5ex]{
    % Define the nodes with 'A' in smaller font
    \node[draw, circle, inner sep=1pt] (tensor1) {\scriptsize A};
    \node[draw, circle, inner sep=1pt, right=0.3cm of tensor1] (tensor2) {\scriptsize A};
    \node[draw, circle, inner sep=1pt, right=0.3cm of tensor2] (tensor3) {\scriptsize A};
    \node[draw, circle, inner sep=1pt, right=0.3cm of tensor3] (tensor4) {\scriptsize A};
    \node[draw, circle, inner sep=1pt, right=0.3cm of tensor4] (tensor5) {\scriptsize A};

    % Draw the edges between tensors
    \draw (tensor1) -- (tensor2);
    \draw (tensor2) -- (tensor3);
    \draw (tensor3) -- (tensor4);
    \draw (tensor4) -- (tensor5);

    % Add vertical legs to each tensor
    \draw (tensor1) -- ++(0,0.3);
    \draw (tensor2) -- ++(0,0.3);
    \draw (tensor3) -- ++(0,0.3);
    \draw (tensor4) -- ++(0,0.3);
    \draw (tensor5) -- ++(0,0.3);
}.  Similarly, an MPO will be represented as: 
\tikz[baseline=-0.5ex]{
    % Define the nodes with 'A' in smaller font
    \node[draw, circle, inner sep=1pt] (tensor1) {\scriptsize M};
    \node[draw, circle, inner sep=1pt, right=0.3cm of tensor1] (tensor2) {\scriptsize M};
    \node[draw, circle, inner sep=1pt, right=0.3cm of tensor2] (tensor3) {\scriptsize M};
    \node[draw, circle, inner sep=1pt, right=0.3cm of tensor3] (tensor4) {\scriptsize M};
    \node[draw, circle, inner sep=1pt, right=0.3cm of tensor4] (tensor5) {\scriptsize M};

    % Draw the edges between tensors
    \draw (tensor1) -- (tensor2);
    \draw (tensor2) -- (tensor3);
    \draw (tensor3) -- (tensor4);
    \draw (tensor4) -- (tensor5);

    % Add vertical legs to each tensor
    \draw (tensor1) -- ++(0,0.3);
    \draw (tensor2) -- ++(0,0.3);
    \draw (tensor3) -- ++(0,0.3);
    \draw (tensor4) -- ++(0,0.3);
    \draw (tensor5) -- ++(0,0.3);
        % Add vertical legs to each tensor
    \draw (tensor1) -- ++(0,-0.3);
    \draw (tensor2) -- ++(0,-0.3);
    \draw (tensor3) -- ++(0,-0.3);
    \draw (tensor4) -- ++(0,-0.3);
    \draw (tensor5) -- ++(0,-0.3);
}. This form of graphical representation is particularly powerful for reasoning about complex tensor networks, as it avoids cumbersome algebraic notation while retaining structural clarity.

\subsection{The quantized tensor train} 

The quantum-inspired approach to solving partial differential equations builds on the quantized tensor train (QTT) formalism. For illustrative purposes, consider a one dimensional function defined on the unit interval, $f:[0,1[\rightarrow \mathbb{R}$.  The unit interval is discretized  into $N=2^n$ equally distributed points. We will typically omit the $x=1$ point. Any grid value $x$ can be conveniently written in its binary expansion as

\begin{equation}
    x= 0.x_1x_2\hdots x_c = \frac{x_1}{2} + \frac{x_2}{2^2} + \cdots +\frac{x_c}{2^c},
\end{equation}
in terms of bits $x_i=\{0,1\}$. Similarly the function $f$ can be represented as a tensor 
\begin{equation}
f_{x_1x_2...x_c} = f(0.x_1x_2...x_c),
\end{equation}
with $n$ indices, each taking values $x_i=\{0,1\}$. The main idea of the QTT construction is to approximate $f$ as a tensor train. Because the indices $i=1,...,c$ reflect length scales of the unit interval, the QTT representation explicitly encodes the multiresolution nature of the function. The virtual dimension of the tensor train reflects correlations between the coarse and fine degrees of freedom in the system. As a consequence, many piecewise smooth functions can be represented efficiently as a QTT \cite{kazeev2012}. 

If a function can be represented faithfully as a QTT with constant bond dimension $R$, then the memory cost of the function on the equidistant grid is $\mathcal{O}(R^2 c)$, as opposed to $N=2^c$ for a dense representation.
There are multiple ways of extending the QTT to higher dimensions. 
As an illustrative example, consider the tensor $\T^{16  \times 16 \times 16}$ representing the discretization of ${f(x,y,z)=
\sin(\alpha_1 x + \phi_1)
\sin(\alpha_2 y + \phi_2)
\sin(\alpha_3 z + \phi_3)}$, where $\phi$ is an arbitrary phase on $(0,1)^3$ with $16$ points in each dimension. Using the QTT decomposition, this tensor can be exactly expressed as rank-2 tensor of the form:

\tikz[baseline=-0.5ex]{
    % Define the nodes with 'A' in smaller font
    \node[draw, circle, inner sep=1pt] (tensor1) {\phantom{\scriptsize Q}};
    \node[draw, circle, inner sep=1pt, right=0.3cm of tensor1] (tensor2) {\phantom{\scriptsize Q}};
    \node[draw, circle, inner sep=1pt, right=0.3cm of tensor2] (tensor3) {\phantom{\scriptsize Q}};
    \node[draw, circle, inner sep=1pt, right=0.3cm of tensor3] (tensor4) {\phantom{\scriptsize Q}};
    \node[draw, circle, inner sep=1pt, right=0.3cm of tensor4] (tensor5) {\phantom{\scriptsize Q}};
    \node[draw, circle, inner sep=1pt, right=0.3cm of tensor5] (tensor6) {\phantom{\scriptsize Q}};
    \node[draw, circle, inner sep=1pt, right=0.3cm of tensor6] (tensor7) {\phantom{\scriptsize Q}};
    \node[draw, circle, inner sep=1pt, right=0.3cm of tensor7] (tensor8) {\phantom{\scriptsize Q}};
    \node[draw, circle, inner sep=1pt, right=0.3cm of tensor8] (tensor9) {\phantom{\scriptsize Q}};
    \node[draw, circle, inner sep=1pt, right=0.3cm of tensor9] (tensor10) {\phantom{\scriptsize Q}};
    \node[draw, circle, inner sep=1pt, right=0.3cm of tensor10] (tensor11) {\phantom{\scriptsize Q}};
    \node[draw, circle, inner sep=1pt, right=0.3cm of tensor11] (tensor12) {\phantom{\scriptsize Q}};

     % Draw the edges between tensors
    \foreach \i [evaluate=\i as \j using int(\i+1)] in {1,...,11} {
        \draw (tensor\i) -- (tensor\j);
        \path (tensor\i) -- (tensor\j) coordinate[midway] (mid\i);
        %\node[above=0.2pt] at (mid\i) {\scriptsize 2};
    }

    % Draw vertical legs and annotate each with "2"
    \foreach \i in {1,...,12} {
        \draw (tensor\i) -- ++(0,0.5);
    }
    \foreach \i in {1,...,4} {
        \node[above] at ([yshift=0.4cm]tensor\i) {\scriptsize \textcolor{red}{2}};
    }
        \foreach \i in {5,...,8} {
        \node[above] at ([yshift=0.4cm]tensor\i) {\scriptsize \textcolor{nicegreen}{2}};
    }
        \foreach \i in {9,...,12} {
        \node[above] at ([yshift=0.4cm]tensor\i) {\scriptsize \textcolor{blue}{2}};
    }
}
\\~ 
\noindent where we only need $48$ evaluations of the function, instead of $4096$, to build this tensor. This construction is given in Appendix \ref{appen:analytic_qtt_sine} and we place each dimension in a \textit{serial} ordering.

Besides the discrete sine function other vectors and matrices allow for explicit low-rank QTT representations. An especially important case for us is the  tridiagonal matrix %band diagonal matrix 
\cite{kazeev2012}: 
\begin{lemX}\label{lemma_1}
Let $I = 
\begin{psmallmatrix}
1 & 0\\
0 & 1
\end{psmallmatrix}$, 
$J = 
\begin{psmallmatrix}
0 & 1\\
0 & 0
\end{psmallmatrix}$, 
$J' = 
\begin{psmallmatrix}
0 & 0\\
1 & 0
\end{psmallmatrix}$, and $\alpha, \beta, \gamma\in \mathbb{C}$, then for any integer 
$c \geq 2$, the $2^c \times 2^c$ matrix 
\begin{align*}
D_{\alpha,\beta,\gamma} = 
\begin{pmatrix}
\alpha & \beta &  &  &    \\
\gamma & \alpha & \beta &  &    \\
 & \ddots & \ddots & \ddots &    \\
\end{pmatrix}
\end{align*}
of size  has an explicit QTT representation with bond dimension $3$, given by:
\begin{align*}
    D_{\alpha,\beta,\gamma} \! = \!
    \begin{bNiceMatrix}
        \I & \J' & \J
    \end{bNiceMatrix} \! \bowtie \!
    \begin{bNiceMatrix}
        \I & \J' & \J\\
          & \J  &   \\
          &    & \J'  \\
    \end{bNiceMatrix}^{\bowtie (c-2)} \!\!\!\!\bowtie\!
    \begin{bNiceMatrix}
        \alpha \I + \beta \J + \gamma \J'\\
          \gamma \J  \\
          \beta \J' \\
    \end{bNiceMatrix}.
\end{align*}
\end{lemX}
The \textit{inner core product} denoted by $\bowtie$ is defined as:
\begin{align*}
\T \bowtie \G &= 
\begin{bNiceMatrix}
\T_{11} & \T_{12}\\
\T_{21} & \T_{22}\\
\end{bNiceMatrix} \bowtie
\begin{bNiceMatrix}
\G_{11}\\
\G_{21}\\
\end{bNiceMatrix} \\ &=
\begin{bNiceMatrix}
T_{11} \otimes \G_{11} +  T_{12} \otimes \G_{21}\\
T_{21} \otimes \G_{11} +  T_{22} \otimes \G_{21}
\end{bNiceMatrix},
\end{align*}
for  $\T$ a $(2,2,2,2)$ tensor  and $\G$ a $(2,2,2,1)$ tensor.

\subsection{Low-rank QTT Representation of a Function}\label{sec:qtt_rep_fun}

When numerically solving PDEs, we require an efficient way to encode bounded continuous functions in the QTT format. More specifically, to match the construction of differential elements presented in the next section, we seek a QTT representation such that the discretization of $f(x)$ in the interval $(0,1)$ with $2^c$ points is given by a QTT with $c$ cores, preferably maintaining a low-rank structure.

There exist multiple approaches for constructing such a representation. Some elementary functions, such as exponentials  and polynomials, admit explicit analytic low-rank QTT representations. Additionally, various numerical techniques can be employed, including the  TT-SVD algorithm \cite{tt_ose},  tensor cross interpolation (TCI) \cite{oseledets2010tt,dolgov2020parallel,QTCI} and  the multiscale interpolative construction  \cite{lindsey2023multiscale}.

Our framework supports these approaches, each of which has its own advantages and limitations. The choice of method depends on the structure of the function and also on whether one is aiming for accuracy or speed.

\paragraph{Explicit Analytic Low-Rank QTT Representations}
Certain elementary functions have well-established analytic QTT constructions that result in low-rank tensor representations. For instance, it is known that the exponential function $f(x) = e^{\alpha x}$ and the sine function $f(x) = \sin(\alpha x + \phi)$ (given in Appendix \ref{appen:analytic_qtt_sine}) admit explicit QTT representations of ranks $1$ and $2$, respectively, independent of the grid resolution. Polynomial functions admit a QTT representation with a rank equal to the degree of the polynomial plus one. Some of these constructions can be found in Section 2.2 of \cite{dolgov_phd}.

\paragraph{TT-SVD}
The TT-SVD and its variants are widely used methods for obtaining a low-rank QTT representations of functions. The process involves first forming the full tensor representation of the discretized function and then applying a sequential SVD with truncation at each step to obtain a low-rank QTT format. While TT-SVD guarantees an optimal rank reduction in terms of Frobenius norm error, it requires constructing and storing the full tensor, making it computationally expensive for high-dimensional problems.

\paragraph{Tensor cross interpolation} The TCI algorithm generalizes matrix cross approximation to tensors by implicitly constructing a QTT representation through oracular access to the function. It selects a small, structured subset of tensor entries, enabling efficient approximation when the full tensor is too large to store or evaluate.
This approach assumes the ability to query the function at arbitrary multi-indices—typically feasible in synthetic PDE problems with known models, but often unavailable in data-driven settings limited to fixed or noisy samples. Consequently, TCI is less suited for learning from data but highly effective when such access exists.
In such cases, TCI offers substantial computational and memory savings and has proven effective in applications like high-dimensional PDEs, uncertainty quantification, and surrogate modeling. 

\paragraph{Multiscale Interpolative QTT Construction}\label{sec:interpolation}
The multiscale interpolative construction of QTT \cite{lindsey2023multiscale} provides an alternative explicit construction based on interpolation. In essence, this method constructs tensor cores by evaluating the target function at $M$ interpolation nodes on a grid (in our case, a Chebyshev-Lobatto grid). 
The resulting QTT has ranks of size $M+1$. As opposed to the TCI implicit algorithm, the interpolative algorithm allows for better control of the error, and is more resilient to noisy data. 
In Appendix~\ref{appen:heat_2d_td_bc}, we demonstrate an additional advantage of our QTT framework equipped with this method by solving the 2D heat Equation with complex time-dependent boundary conditions using this interpolative QTT construction.

Each of these methods provides a viable approach to constructing QTT representations, and the optimal choice depends on the specific function properties and computational goals and constraints in a given application.

\subsection{Finite Difference in QTT format}\label{sec:fin_dif_qtt}

We review the finite-difference discretization in QTT format. For simplicity, we will consider a linear homogeneous second order PDE in one variable:

\begin{align*}
      p\frac{\partial^2}{\partial x^2}u(x) + s\frac{\partial}{\partial x}u(x) + v u(x) = f(x),
\end{align*}
with $x \in (a,b)$, constants $p,s$ and $v$ with Dirichlet boundary conditions $u(a) = 0$ and $u(b) = 0$.
This equation can be discretized into $2^c$ points
giving
\begin{align*}
    &p \left(U_{i-1} - 2U_i + U_{i+1}\right)/h^2\\&+ s \left( U_{i+1} - U_{i-1}\right)/2h  + vU_i = f_i,
\end{align*}
where $h = (b-a)/2^c$ and, for example, $U_i \approx u(x_i)$ with $a \leq i \leq b$. We can write this equation as the linear system $M^{(c)}_{p,s,v}U=F$ where 
\begin{align*}
M_{p,s,v} = 
\begin{pmatrix}
\alpha & \beta &  &  &    \\
\gamma & \alpha & \beta &  &    \\
 & \ddots & \ddots & \ddots &    \\
\end{pmatrix},
\quad
\begin{aligned}
&\alpha = h^2v - 2p, \\
&\beta = p + hs/2, \\
&\gamma = p - hs/2,
\end{aligned}
\end{align*}
is a $2^c \times 2^c$ tridiagonal matrix, $U = \begin{bmatrix} U_1& U_2 & \cdots &U_{2^c} \end{bmatrix}^{T}$ and $F = h^{-2} \begin{bmatrix}f_1 & f_2 & \cdots & f_{2^c} \end{bmatrix}^{T}$.
By Lemma \eqref{lemma_1}, the matrix $M_{p,s,v}$ has an explicit QTT representation of maximal rank $3$. We can build a $c$ cores QTT representation of the vector $F$ using one of the techniques presented in Section \ref{sec:qtt_rep_fun}.

By means of tensor product we can generalize this 1D scheme to higher dimensions. For the 2D case consider the second-order partial differential equation 
\begin{align*}
    p u_{xx} + q u_{yy} + r u_{xy} + 
    s u_x + t u_y + v u = f(x,y).
\end{align*}
Doing a similar discretization of this PDE into $2^{2c}$ points we can write it as a linear system $M^{(2c)}_{p,q,r,s,t,v} U_2 = F_2$, with
\begin{align*}
    M^{(2c)}_{p,q,r,s,t,v} &= (M^{(c)}_{p,s,v} \otimes \I_{2^c}) + (\I_{2^c} \otimes M^{(c)}_{q,t,0})\\ &+ 
    \left( (M^{(c)}_{0,r,0} \otimes \I_{2^c}) \times (\I_{2^c} \otimes M^{(c)}_{0,1,0}) \right),
\end{align*}
where $\I_{2^c}$ is the $2^c \times 2^c$ identity matrix and $U_2$ and $F_2$ are defined accordingly. A key observation is that we can do this same construction for PDEs in higher dimensions and always get a low-rank QTT representation of the matrix ``$M$'' and vector $F$ independently of the number of discretization points.

\subsection{Solving a system of linear equations - ALS}\label{sec:als}
Now that we have presented how to discretize the PDE and construct the QTT representation of its elements, we need an efficient way to solve the resulting system of linear equations. In this section, we provide a brief overview of the Alternating Linear Scheme (ALS) \cite{als_main}, which is a fundamental method for solving linear systems in the TT format.
The ALS method is an iterative optimization approach where tensor cores are updated sequentially through alternating sweeps. Each iteration involves contracting the tensor networks to obtain a local system of equations, solving for an optimal update of a single core while keeping the others fixed. Given an initial guess for the solution, ALS updates the TT cores successively in a bidirectional manner (left-to-right and then right-to-left), corresponding to one full sweep. MALS (modified alternating linear scheme) works similarly but optimizes two cores at a time; in practice, this usually implies higher accuracy but with a small increase in overall run time.
One of the main advantage of (M)ALS is that is extremely fast compared to other optimization methods, while naturally incorporating pre-conditioning, by restricting the learning to the manifold of QTT functions with bounded bond dimension. 

To solve a linear system of the form $Ax = b$, we use the following function: $\text{(M)ALS}(A_{\text{QTT}}, \hat{x}_{\text{QTT}}, b_{\text{QTT}}, \textit{sweeps})$, where $A_{\text{QTT}}$ and $b_{\text{QTT}}$ are the QTT representation of the matrix $A$ and vector $b$ respectively and \textit{sweeps} the number of full sweeps. The term $\hat{x}_{\text{QTT}}$ represents an initial guess for the solution in QTT format. The ranks of the initial guess $\hat{x}_{\text{QTT}}$ plays a crucial role in the efficiency and accuracy of the (M)ALS method. A low-rank initial guess can significantly speed up convergence, but if the ranks are too low, the algorithm may fail to achieve the desired accuracy. Conversely, overestimating the ranks increases computational cost unnecessarily.

There are multiple ways to construct $\hat{x}_{\text{QTT}}$. A simple approach is to use the same rank structure of $b_{\text{QTT}}$, possibly adding a constant factor to its ranks. Alternatively, one could start with a random QTT such that the ranks increase by following an arithmetic (or geometric) progression until the middle core and then decrease symmetrically toward the final core. Our framework implements all these strategies, allowing for easy adjustments to balance efficiency and accuracy.
In \cite{mgr}, the authors adopt a DMRG-like approach to solving the system of linear equations. Their method initializes $\hat{x}_{\text{QTT}}$ using a coarser grid solution, which is then mapped to a finer grid via a prolongation MPO. However, in our experiments, this strategy did not lead to significant improvements in either accuracy or runtime for our method.

The computational complexity of ALS is $\mathcal{O}(c \gamma r^3 R^2 n^2)$, where $c$ is the number of cores in $A_{\text{QTT}}$, $\gamma$ is the number of iterations required to solve the local system of equations, $r$ is the maximum rank of either $\hat{x}_{\text{QTT}}$ or $b_{\text{QTT}}$, $R$ is the maximum rank of $A_{\text{QTT}}$, and $n$ is the maximum mode size (in our framework $n = 2$). In contrast, the dominant computational complexity of the best classical methods is of order $\mathcal{O}(2^c)$. As a final remark we note that, since $A_{\text{QTT}}$ admits an exact low-rank representation and we often also obtain a low-rank representation for $b_{\text{QTT}}$, the overall complexity can remain polynomial in $c$, enabling exponential speedup over classical methods. This allows us to handle arbitrarily fine discretizations, in principle.

\section{Finite Difference Method for the 2D Poisson in QTT Format}\label{sec:fd_poisson}
In this section, we focus on solving the 2D Poisson equation using a finite difference scheme within the QTT framework. This equation serves as a fundamental building block for addressing more complex PDEs, including nonlinear, time-dependent, or higher-dimensional cases.

Consider the 2D Poisson equation
\begin{align*}
      &\Delta u = \frac{\partial^2 u}{\partial x^2} + \frac{\partial^2 u}{\partial y^2} =  f(x,y), 
\end{align*}
where $(x,y) \in \Omega = (a,b) \times (d,e)$, subject to Dirichlet boundary conditions $u(x,y) = g(x,y)$ if $(x,y) \in \partial\Omega$.
We start by discretizing the region $\Omega$ into $\Omega_h$ with $N = 2^c$ points in each spatial dimension. Each grid point in the $x-$dimension is given by $x_i = a + ih$, for ${i=0,\hdots,N+2}$, with $h_x = (b-a)/(N+2)$, similarly we define $y_j$. Let $w_{ij} := u(x_i,y_j)$ represent the approximate solution at each grid point. Similarly, $f_{ij}:= f(x_i,y_j)$. Using a second-order central difference approximation, we can discretize this PDE wlog over a square region:
\begin{align}\label{discrete_poisson}
\begin{split}
    \frac{1}{h^2} \big[&w_{i+1, j} - 2 w_{i, j} + w_{i-1, j}\\ &+ w_{i, j+1} - 2 w_{i, j} + w_{i, j-1}\big] = f_{ij},
    \end{split}
\end{align}
for $i,j = 1, \hdots, N$. Equation \eqref{discrete_poisson} can be written as a system of linear equations
\begin{align}\label{linear_system_1}
    A w = -h^2f - b,
\end{align}
where  $A = \Delta_{DD} \otimes \I_{2^c} + \I_{2^c} \otimes \Delta_{DD}$, with $\Delta_{DD} = D^{(c)}_{-2,1,1}$. This matrix is the 2D discrete Laplacian and by Lemma \eqref{lemma_1} has a exact low-rank QTT representation with bond dimension $6$. In the serial ordering of the tensors, the vectors $w$ and $f$ can be written as:
\begin{align}\label{w_2d}
    w &= \begin{bmatrix} w_{1,1} & \cdots & w_{N,1} &  \cdots\cdots & w_{N,2} & \cdots  & w_{N,N} \end{bmatrix}^T,
\end{align}
\begin{align}\label{f_2d}
    f &= \begin{bmatrix} f_{1,1} & \cdots & f_{N,1} &  \cdots\cdots & f_{N,2} & \cdots &  f_{N,N} \end{bmatrix}^T.
\end{align}
The $b$ vector accounts for the boundary conditions. One way to define $b$ is as follows:
\begin{align}\label{b_vec}
\begin{split}
b = 
&(b^{(\text{left})} \otimes | 0 \rangle) + 
(b^{(\text{right})} \otimes | 1 \rangle)\\   
&+ (| 0 \rangle \otimes b^{(\text{bottom})}) + 
(| 1 \rangle \otimes b^{(\text{top})}),
\end{split}
\end{align}
where 
\begin{align*}
    b^{(\text{left})} &= \begin{bmatrix}
        w_{0,1} \\ w_{0,2} \\ \vdots \\ w_{0,N}
    \end{bmatrix}, \quad     
    &b^{(\text{right})} = \begin{bmatrix}
        w_{N,1} \\ w_{N,2} \\ \vdots \\ w_{N,N}
    \end{bmatrix}, \\
    b^{(\text{bottom})} &= \begin{bmatrix}
        w_{1,0} \\ w_{2,0} \\ \vdots \\ w_{N,0}
    \end{bmatrix}, \quad     
    &b^{(\text{top})} = \begin{bmatrix}
        w_{1,N} \\ w_{2,N} \\ \vdots \\ w_{N,N}
    \end{bmatrix},
\end{align*}
here, $| 0 \rangle$ and $| 1 \rangle$ are the vectors of length $N$ of the form $(1,0,\hdots,0)$ and $(0,0,\hdots,1)$, respectively. To build the QTT representation of the four $b^{(\cdot)}$ vectors and $f$ we use one of the methods presented in \ref{sec:qtt_rep_fun}. The QTT representation of $| 0 \rangle$ and $| 1 \rangle$ is given in Appendix \ref{appen:1d_bc_qtt}.

With all required components represented in the QTT format, we are now in position to use (M)ALS to solve the system of linear equations \eqref{linear_system_1}.

\subsection{Scaling Comparison: QTT Solver vs. Algebraic Multigrid}\label{sec:Laplace_Poisson}
In this subsection, we compare the performance of our QTT-based solver with the widely used Algebraic Multigrid (AMG) method, implemented by the PyAMG library. AMG is an industry standard for solving PDEs, known for its linear complexity, $\mathcal{O}(N)$, where
$N$ is the number of discretization points.  However, for high-dimensional or high resolution PDEs, AMG succumbs to the curse of dimensionality as $N$ grows exponentially. In contrast, our QTT solver demonstrates a scaling behavior of $\mathcal{O}(\log(N))$. 
This remarkable speedup has its origin in three key features of the QTT method: (i) there exists a low-rank QTT representation of discrete differential elements (ii) that the QTT can represent piecewise smooth functions with low bond dimension, (iii) that the (M)ALS solver acts as an effective preconditioner.

Our results highlight the competitiveness of the QTT solver, particularly in cases where AMG's efficiency diminishes. Furthermore, our framework is highly flexible, allowing users to easily adjust the trade-off between time and accuracy by tuning various parameters of the solver. For completeness, additional benchmarks for the 2D Poisson equation are provided in the appendix.

We now turn to the specific example of the 2D Laplace equation: $\Delta u=0$ in the domain $(0,1) \times (0,1)$ with all zero Dirichlet boundary conditions except for ${u(x,0) = \sin(k \pi x) \sinh(k \pi)}$. The corresponding analytical solution is
$u(x,y) = \sin(k \pi x) \sinh(k\pi(1-y))$. For values of $k$ larger than $2$, the solution shows a sharp kink at the origin (see Appendix \ref{appen:problem_1}), requiring a very fine grid to obtain a high precision solution.  

\begin{figure}[ht]
    \centering
    \includegraphics[width=\linewidth]{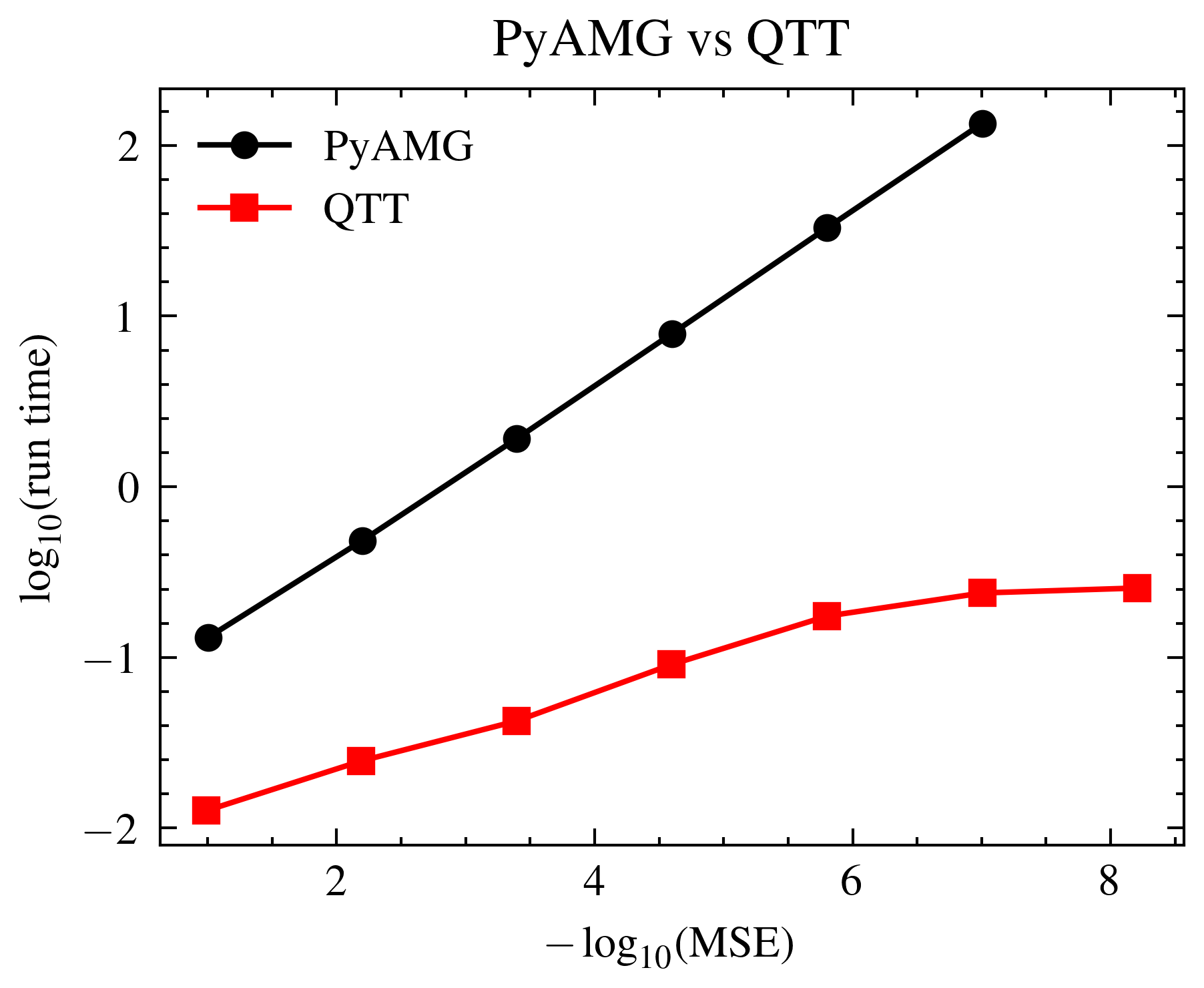}
    \caption{Log-log plot comparing the run time and accuracy between PyAMG and QTT when a low rank representation of the boundary term is available. }
    \label{fig:log-log Poisson 2d}
\end{figure}

In Figure \ref{fig:log-log Poisson 2d}, we compare the performance of one of our QTT-Solvers against PyAMG for $k=3$. The  QTT  solver used an exact analytical low-rank QTT representation of the non-zero boundary condition. PyAMG was configured to use a standard Ruge-Stüben algebraic multigrid method. The complete data for this plot are provided in the Benchmark Section \ref{app:Benchmarks}. The small variations in runtime and MSE observed in our QTT solver are attributed to the inherent randomness in the initial guess for the solution.

As expected, PyAMG has a runtime that scales linearly as we increase the number of discretization points, whereas our QTT solver scales logarithmically with the number of discretization points. We note that to achieve an MSE of $10^{-6}$, PyAMG took approximately $30$ seconds, while our QTT solver took only about $0.2$ seconds. The last red node indicates that the QTT solver requires just another $0.1$ seconds to reach an MSE of $10^{-10}$.
There are several strategies to improve the solution accuracy. As a rule of thumb, updating two cores at a time (MALS) tends to be slower than updating one core at a time (ALS), but often results in higher precision or at least matches with ALS. Increasing the bond dimension or choosing a more appropriate initial guess of the solution could also result in higher accuracy. We illustrate one of these strategies on Problem \ref{appen:problem_3} in the benchmark section and also introduce a standard setup that performs well across various elliptic PDEs.

To further highlight the advantages of our QTT-based methods for high-dimensional PDEs, we now solve the anisotropic 3D Poisson equation:
\begin{align*}
    \frac{\partial^2u}{\partial x^2} + \epsilon_1 \frac{\partial^2u}{\partial y^2} + \epsilon_2 \frac{\partial^2u}{\partial z^2}&= -\sin(\pi x) \sin(\pi y) \sin(\pi z),
\end{align*}
with $(x,y,z) \in (0,1)^3$ and all zero boundary conditions. The analytical solution is given by:
\begin{align*}
    u(x,y) = \dfrac{\sin(\pi x) \sin(\pi y) \sin(\pi z)}{\pi^2 (1 +\epsilon_1 + \epsilon_2)}.
\end{align*}

In the following table we report the runtime comparing PyAMG and our QTT solver to achieve a certain level of accuracy. For completeness, in Benchmark Section \ref{appen:problem2}, we present results for the case $\epsilon_1 = \epsilon_2 = 1$, where both solvers reach the discretization error while exhibiting the same runtime behavior as discussed below. 

\begin{table}[h]
    \centering
    \begin{tabular}{cccc}
        \toprule
         \multicolumn{2}{c}{\textbf{PyAMG}} & \multicolumn{2}{c}{\textbf{QTT}} \\ 
         \cmidrule(lr){1-2} \cmidrule(lr){3-4} 
        Time(s)   & MSE        & Time(s)   & MSE \\ 
        \hline
        \midrule
        0.011     & 2.81e-06   & 0.001     & 2.81e-06  \\
        0.021     & 1.90e-07   & 0.001     & 1.90e-07  \\
        0.056     & 1.24e-08   & 0.002     & 1.24e-08  \\
        0.32     & 8.10e-10   & 0.005     & 8.01e-10   \\
        2.62     & 5.08e-11   & 0.010     & 5.08e-11   \\
        22.39    & 3.20e-12   & 0.020     & 3.20e-12   \\
        347.34   & 2.00e-13   & 0.031     & 2.00e-13   \\
        -         & -          & 0.053     & 1.25e-14  \\
        \hline
        \bottomrule
    \end{tabular}
    \caption{Performance Comparison of PyAMG and our QTT Poisson solver with $\epsilon_1 = 0.001$ and $\epsilon_2 = 0.0001$ from $2$ to $(8)$ $9$ cores per dimension.}
\end{table}\label{tab:discretization}

The presented QTT method get more accurate as we increase the number of discretization points while maintaining logarithmic runtime scaling.  
These experiments confirm that the QTT method can exponentially outperform state of the art classical linear solvers, when the source or boundary terms have explicit low rank constructions. When this is not the case, the dominant error term can come from the QTT construction of the source terms. In Section \ref{sec:data}, we explore the case when the source terms are learned from data.

\section{Finite Difference Methods for Burgers' Equation in QTT Format}\label{sec:fd_burgers}

Most PDEs of practical relevance are time-dependent and nonlinear. In this section, we study Burgers equation, and show that the exponential advantage seen for the Poisson equation extends to time-dependent nonlinear problems. 
We start by building a time-stepping finite difference scheme within the QTT framework to solve Burgers' equation, aiming to show how effectively our approach handles nonlinearity and different boundary conditions. Next, leveraging the low-rank structure of the QTT format, we extend this method into a space-time QTT solver by treating the time dimension as a spatial one. We end this section by comparing the runtime and accuracy of these two methods on a specific example of Burgers' equation. See also the space-time analysis of the heat equation in Appendix \ref{appen:heat_ts_st}.

It is important to pause at this point and discuss why nonlinearity and time dependence pose significant challenges. Nonlinear, time-dependent PDEs are particularly demanding for numerical algorithms due to the combined difficulties of maintaining stability, ensuring solver convergence, and accurately capturing evolving dynamics. Nonlinearities can make standard schemes unstable or require iterative solvers that may fail to converge without good initial guesses. Simultaneously, time dependence requires careful time-stepping to resolve transient behavior and prevent error accumulation over long simulations. These factors make such equations particularly demanding for numerical algorithms.

\subsection{Time-stepping QTT Solver}
In this subsection, we present a time-stepping QTT solver applied to the 1D + 1t Burgers' equation with Dirichlet boundary conditions on both ends, followed by a case with Neumann-Dirichlet boundary conditions. This method is easily extendable to higher spatial dimensions and other boundary conditions. 

This equation is given by:
\begin{align}\label{burgers_eq}
      &\frac{\partial u}{\partial t} =  \nu \frac{\partial^2 u}{\partial x^2} - u\frac{\partial u}{\partial x}, 
\end{align}
with $(x,t) \in \Omega = \{a \leq x \leq b, t \geq 0\}$, initial condition $u(x,0) = g_0(x)$ , and for now, we assume Dirichlet boundary conditions: $u(a,t)=g_1(t)$ and $u(b,t)=g_2(t)$. 

Let $w_{ij} \approx u(x_i,t_j)$ and consider the following discretization scheme:

\begin{align*}
    (w_{i,j+1} - w_{i,j})/l =  
    &\nu \underbrace{(w_{i-1,j+1} - 2w_{i,j+1} + w_{i+1,j+1})/h^2}_{\alpha}\\
    &-u \underbrace{(w_{i-1,j+1} - w_{i+1,j+1})/(2h)}_{\beta},
\end{align*}
where $l$ is the time step size and $h = (b-a)/N$, with $N=2^c$, the dimension step size. More compactly,
\begin{align}\label{burgers_1_disc}
   -\nu \alpha l + w_{i,j+1} + u \beta l &= w_{i,j}.
\end{align}
Assuming all time steps have the same size we will approximate $u$ by a linear interpolation ${u_{i,j+1} = 2u_{i,j} - u_{i,j-1}}$. Substituting in Equation \eqref{burgers_1_disc} we have
\begin{align}\label{burgers_2_disc}
   -\nu \alpha l + w_{i,j+1} + (2w_{i,j} - w_{i,j-1}) \beta l &= w_{i,j}.
\end{align}
Our goal is to write \eqref{burgers_2_disc} as a system of linear equations to be solved for each time step stating from the solution at the previous one. Typically, this involves treating the nonlinear term at iteration $j+1$ as a constant multiplicative factor.  Specifically, we define the diagonal matrix ${D_j = \textrm{diag}(w_{1,j}, w_{2,j}, \hdots, w_{N,j})}$ and for simplicity of notation, let $D'_j = 2D_j - D_{j-1}$. Now we can approximate Equation \eqref{burgers_2_disc} as the linear system ${(A+ D'_jB) w_{j+1}= w_j + b_{j+1}}$, where
\begin{align*}
A =
\begin{bNiceMatrix}
1-2r & r &  &  &    \\
r & 1-2r & r &  &    \\
& \ddots & \ddots & \ddots     \\
\end{bNiceMatrix}_{(N \times N)},
\end{align*}

\begin{align*}
w_j =
\begin{bNiceMatrix}
    w_{1,j} \\ w_{2,j}\\ \Vdots \\ w_{N,j}
\end{bNiceMatrix},~
    b_{j+1} = 
\begin{bNiceMatrix}
    r w_{0,j+1} \\ 0 \\ \Vdots \\ 0 \\ r w_{N,j+1}
\end{bNiceMatrix},
\end{align*}
with $r=-\nu l/h^2$ and
\begin{align*}
B = \frac{l}{2h}
\begin{bNiceMatrix}
 0 & 1    &0   &  &\\
 -1    & 0 & 1  & &\\
  &\Ddots & \Ddots & \Ddots\\
  & &\phantom{-1} &\phantom{0}&\phantom{1}
\end{bNiceMatrix}_{(N \times N)}.
\end{align*}
In this formulation, only $w$ needs to be solved, while $D'$ is treated as a constant. 
Note that we can build the low-rank QTT representation of $B$ using Lemma \eqref{lemma_1} and the QTT representation of $D$ using construction \ref{appen:vec_diag_mpo}.

With a small modification to this method we could also solve Burgers' equation with other types of boundary conditions. For instance, let us consider Neumann-Dirichlet boundary condition $u_x(a,t)=g_1(t)$ and $u(b,t)=g_2(t)$. The only difference from the previous method is that now we will solve the linear system $(A' + D'_jB')w_{j+1} = w_j + b_{j+1}$, where

\begin{align*}
A' = A
-
\begin{bNiceMatrix}
  1 - 2r - \dfrac{1}{h} & r + \dfrac{1}{h} &        &  \\
                        & \ddots           & \ddots &  \\
                        &                  & r      & 1 - 2r - 1
\end{bNiceMatrix},
\end{align*} 
and
\begin{align*}
B' = B -
\begin{bNiceMatrix}
  0   & l/2h   &        &        \\
      & \ddots & \ddots &   \\
      &        & -l/2h & 0
\end{bNiceMatrix}
\end{align*} 
Note that the rightmost matrix on both equations has all its elements equal to zero except for the four entries shown. To build the QTT representation of these matrix we use the construction given in \ref{appen:eraser_qtt}.

The main steps of the time-stepping algorithm to solve Burgers' equation are given below:
\begin{algorithm}[H]
\caption{Time-stepping QTT Solver Burgers' Eq}\label{alg:burgers_time_stepping}
\begin{algorithmic}[1]
\Statex{\textbf{INPUT:}} $c$, $l$, \textit{timesteps}

\Statex{\textbf{OUTPUT:}} $c$ cores QTT representation of the solution at each time step

\Comment{Start first time step}

\Statex{Build the $c$ cores QTT representation of:}

\State $A$ and $B$ \Comment{Lemma \eqref{lemma_1}}

\State $w_0$ from initial conditions \Comment{Any method from \ref{sec:qtt_rep_fun}}

\State $D_0$ using step (2) \Comment{Construction \ref{appen:vec_diag_mpo}}

\State $b_{1}$ \Comment{Construction \ref{appen:1d_bc_qtt}}

\State $\textrm{LHS}_{QTT} \leftarrow A_{\text{QTT}} + D_0 B_{\text{QTT}}$

\State $w_{sol}[0] \leftarrow w_{0_{\text{QTT}}} + b_{1_{\text{QTT}}}$

\State $w_{sol}[1] \leftarrow$ (M)ALS$(\textrm{LHS}_{QTT}, w_{sol}[0], w_{sol}[0])$

\Comment{End first run}

\For{$k=2$ until \textit{timesteps}}

\State Build $D_{k-1}$ using $w_{sol}[k-1]$

\State $\textrm{LHS}_{QTT} \leftarrow A_{\text{QTT}} + (2D_{k-1} + D_{k-2})B_{\text{QTT}}$

\State Build $b_{k_{\text{QTT}}}$ 

\State $\textrm{RHS}_{QTT} \leftarrow w_{sol}[k-1] + b_{k_{\text{QTT}}}$

\State $w_{sol}[k] \leftarrow$ (M)ALS($\textrm{LHS}_{QTT}, \textrm{RHS}_{QTT},\textrm{RHS}_{QTT}$)
\EndFor
\State Return{} $w_{sol}$
\end{algorithmic}
\end{algorithm}

\subsection{Space-time QTT Solver}

As the previous examples illustrate, the low-rank structure of the QTT representation enables us to efficiently handle a high number of discretization points. One possible application of this feature is solving Burgers' equation in the QTT framework by treating the time dimension as a spatial one. By introducing the \textit{runs} parameter, defined as the number of updates of the nonlinear term we can completely eliminate the need for explicit time-stepping, allowing us to solve for the entire solution at once. Through this simple example, we show that it is indeed possible to achieve a logarithmic cost in memory and runtime for time dependent non+linear problems. 

Consider Equation \eqref{burgers_eq} with the same initial and boundary conditions but now treat $t$ as a spatial coordinate, and for simplicity, also in the interval $(a,b)$. Using this formulation, let $w_{ij} \approx u(x_i,t_j)$ and consider the following discretization scheme:
\begin{equation}
\begin{aligned}\label{eq:spacetime_disc}
    \frac{w_{i,j-1} - w_{i,j}}{h} &- \nu
    \frac{- w_{i-1,j} + 2w_{i,j} - w_{i+1,j}}{h^2}\\ 
    &+u\frac{w_{i-1,j} - w_{i,j}}{h}
    = 0,
\end{aligned}
\end{equation}
where wlog $h = (b-a)/N$, with $N=2^c$, is the dimension step size for both $x$ and $t$. 
For simplicity, we assume $a=0$ and $b=1$, which establishes a direct correspondence between the space-time formulation and the time-stepping approach. In theory, under this setup, this means that solving the system in space-time form should yield similar results as performing explicit time-stepping with a time step size of $l=1/N$ for $N$ time steps.

We can rewrite Equation \eqref{eq:spacetime_disc} as a system of linear equations $(A_1 - \nu A_2 + D A_3) w = -b$. Here, the matrices $ A_1, A_2,$ and $A_3$ are tridiagonal Toeplitz matrices, with their construction detailed in Appendix \ref{appen:burgers_st}. The term $D$ is defined similarly to the time-stepping algorithm. The first approximation of $u$ is given by the initial condition $g_0(x)$. Subsequent approximations of $u$ are obtained using the previous approximations of the whole solution. The parameter \textit{runs} specifies the number of iterations used to improve this approximation. The vector $w$ follows the same indexing as the vector defined in Equation \eqref{w_2d}. The right-hand side vector $b$ encodes the discretization of the initial and boundary conditions, which can be efficiently constructed using tensor products and also has a low-rank QTT representation.
Algorithm \ref{alg:burgers_spacetime} presents the main steps of the space-time algorithm to solve Burgers' equation.
\begin{algorithm}[H]
\caption{Space-time QTT Solver  Burgers' Equation}\label{alg:burgers_spacetime}
\begin{algorithmic}[1]
\Statex{\textbf{INPUT:}} $c$, \textit{runs}

\Statex{\textbf{OUTPUT:}} $c$ cores QTT representation of the solution

\Statex{} Build the $c$ cores QTT representation of:

\State $b$ \Comment{Any method from \ref{sec:qtt_rep_fun}}

\State $D_1$ using $b_{\text{QTT}}$ from (1) \Comment{Construction \ref{appen:vec_diag_mpo}}

\State $A_1,A_2$ and $A_3$ \Comment{Lemma \eqref{lemma_1}}

\State $w_{\text{sol}}[0] \leftarrow b_{\text{QTT}}$

\For{$k=1$ until \textit{runs}}

\State $\textrm{LHS}_{QTT} \leftarrow A_{1_{\text{QTT}}} - \nu A_{2_{\text{QTT}}} + D_k A_{3_{\text{QTT}}}$

\State $w_{\text{sol}}[k] \leftarrow$ (M)ALS($\textrm{LHS}_{QTT},w_{\text{sol}}[k-1],b_{QTT}$)

\State Build $D_{k+1}$ using $w_{sol}[k]$

\EndFor
\State Return{} $w_{sol}$
\end{algorithmic}
\end{algorithm}

\subsection{Time Stepping vs Space-time}\label{sec:time_stepping_vs_ST}

To compare both methods, we consider Burgers' equation, given by Equation \eqref{burgers_eq}, with the initial condition  
\begin{align*}
u(x,0) = \frac{2 \nu \pi \sin(\pi x)}{\alpha + \cos(\pi x)},
\end{align*}  
where $\alpha > 1$ and $\nu > 0$, for $x \in (0,1)$. We impose Dirichlet boundary conditions, $u(0,t) = u(1,t) = 0$, for $0 \leq t \leq 1$. The analytical solution is given in \cite{wood_sol_burgers}:  
\begin{align*}
u(x,t) = \frac{2 \nu \pi e^{-\nu \pi^2 t} \sin(\pi x)}{\alpha + e^{-\nu \pi^2 t} \cos(\pi x)}.
\end{align*}
In Appendix \ref{appen:burgers_plots}, we present plots for the space-time method and also the analytical solution for this PDE.

In the table below we compare how long both methods takes to get to approximately the same MSE for different combinations of parameters. For the time-stepping algorithm, we perform $2^7$ time steps of size $1/2^7$, with the spatial dimension discretized into $2^6$ points. In the space-time algorithm, we use a 2D grid with $2^7$ points in each dimension. 
\begin{table}[h]
    \centering
    \begin{tabular*}{0.48\textwidth}{@{\extracolsep{\fill}}cccccccc@{}}
        \toprule
        \multicolumn{2}{c}{\textbf{Parameters}} & \multicolumn{2}{c}{\textbf{Time Stepping}} & \multicolumn{2}{c}{\textbf{Space-Time}} \\ 
        \cmidrule(lr){1-2} \cmidrule(lr){3-4} \cmidrule(lr){5-6}
        $\nu$    & $\alpha$  & Run Time(s) & MSE & Run Time(s) & MSE \\
        \hline
        \midrule
        0.01   & 1.01  & 0.176 & 2.02e-05 & 0.0072 & 3.28e-04 \\
        0.01   & 1.05  & 0.173 & 5.88e-06 & 0.0061 & 6.29e-05 \\
        0.01   & 1.25  & 0.167 & 1.57e-06 & 0.0055 & 5.69e-06 \\
        0.001  & 1.01  & 0.175 & 9.66e-07 & 0.0053 & 1.10e-07 \\
        0.001  & 1.05  & 0.171 & 3.11e-07 & 0.0050 & 4.81e-08 \\
        0.001  & 1.25  & 0.170 & 5.73e-08 & 0.0047 & 2.86e-09 \\
        \hline
        \textbf{1e-05}  & \textbf{1.01}  & - & - & \textbf{0.0141} & \textbf{1.77e-13} \\
        \textbf{1e-07}  & \textbf{1.05}  & - & - & \textbf{0.0116} & \textbf{2.12e-22} \\
        \hline
        \bottomrule
    \end{tabular*}
    \caption{Comparison of time stepping and space-time QTT solvers. In bold we considered $2^{14}$ points in each dimension, demonstrating the the space-time solver can reach regimes unattainable with time-stepping.}
    \label{tab:combined_results}
\end{table}

\begin{figure}[ht]
    \centering
    \includegraphics[width=1\linewidth]{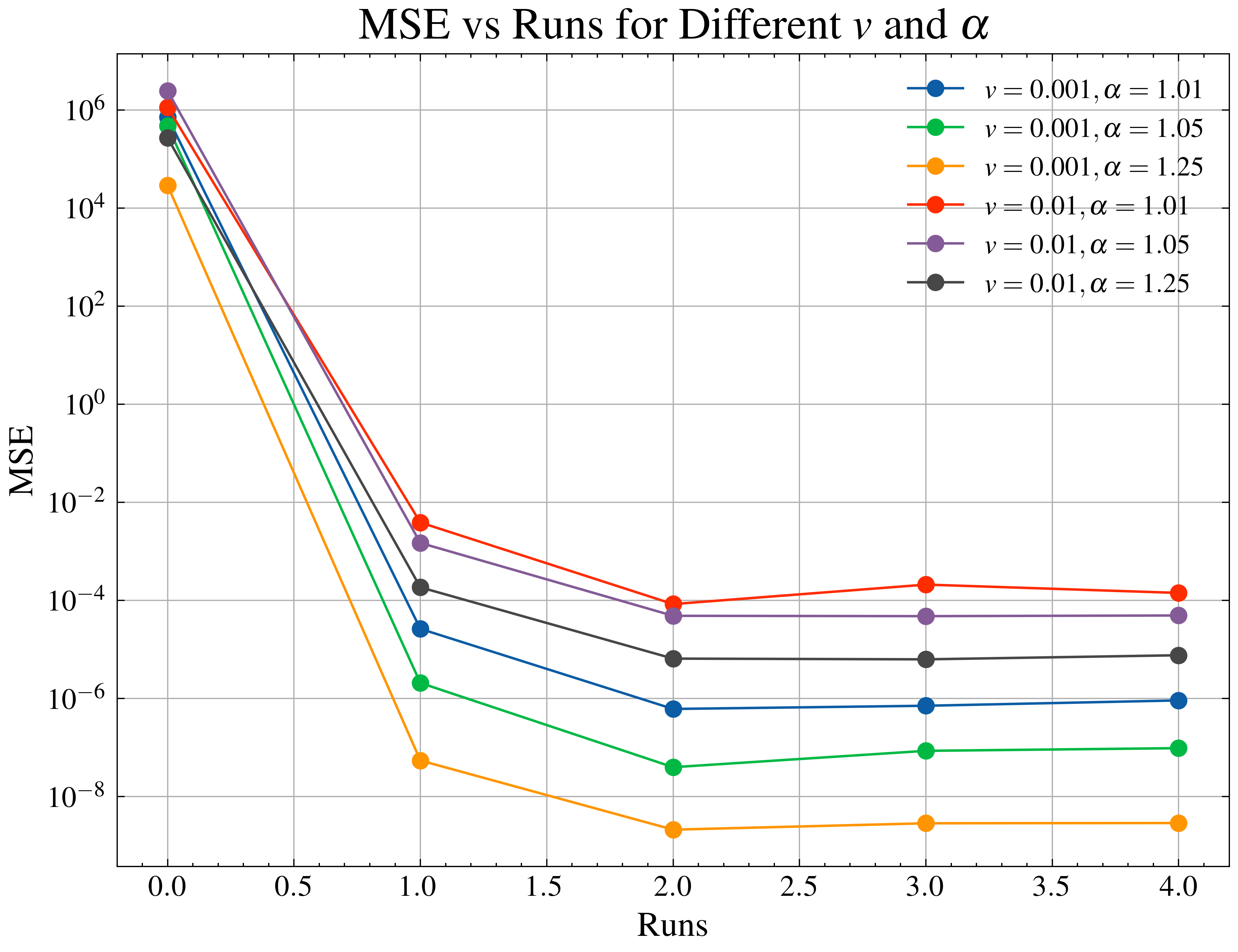}
    \caption{MSE of the space-time QTT solver as a function of the number of \textit{runs}. The \textit{runs} to convergence is essentially independent of the parameters of the PDE.}
    \label{fig:mse_vs_runs}
\end{figure}

For the space-time method, higher accuracy can be achieved by increasing the spatial resolution (as shown in the last two rows of the table). However, unlike the time-stepping approach, the number of \textit{runs} does not need to be increased with the grid resolution. As shown in Figure~\ref{fig:mse_vs_runs}, performing just $2$ \textit{runs} is sufficient to reach the reported accuracy, which remains nearly constant with additional \textit{runs}.

In Appendix \ref{appen:burgers_pinns}, we compare our QTT space-time solver to a well-established benchmark of Burgers’ equation commonly used to evaluate PINNs. Our solver reaches the same level of accuracy while being approximately 100 times faster, highlighting its significant advantage in runtime.

It is important to note here, that the space time solver completely circumvents the CFL condition, meaning that we can reach the $\mathcal{O}(\log(NT))$ total scaling. Finally, we note that using a higher-order discretization scheme, such as Crank-Nicolson, could further improve the accuracy of our QTT solver, particularly for the space-time method.

\newpage
\section{Data-driven approach}\label{sec:data}

Recently, PINNs have emerged as a powerful way to merge machine learning with PDE physics. Encoding PDE residuals, boundary data, and observations in the loss function lets them yield physically consistent solutions from sparse or noisy data. To match the flexibility from PINNs we have to include data points into the QTT framework. Our goal in this section is to highlight how a data-driven approach can be implemented in the QTT framework without having to find a low-rank representation of the function that adequately describes the function passing through all data-points on the boundary condition. We circumvent this problem by introducing splines based on interpolation results from \cite{lindsey2023multiscale}. 

Our method begins by constructing a smooth interpolating function from a given set of data points $\{x_i, y_i\}$. We employ spline interpolation techniques--such as cubic splines or B-splines--to create a continuous function $s(x)$ that accurately approximates the discrete data. The choice between different spline types and degrees allows us to tailor the interpolation to the specific characteristics of the data and the desired smoothness of the function.

Next, we generate a set of interpolation nodes $\{x_j\}$ over the domain of interest, serving as evaluation points for the interpolating function. These nodes can be selected based on various schemes, including Chebyshev nodes, equally spaced nodes, or Legendre nodes, each offering specific advantages. This also introduces a new hyperparameter $M$, the number of interpolation nodes which corresponds to the bond dimension $M+1$ in the interpolation framework. Finally, we use this spline function along with the desired number of cores $c$ for the target QTT and the parameter $M$ as inputs to the interpolation method previously described in Section~\ref{sec:interpolation}. The main steps of this algorithm are outlined below:

\begin{algorithm}[H]
  \caption{Data-Driven QTT Representation via Spline Interpolation}
  \label{alg:qtt_spline_interpolation}
  \begin{algorithmic}[1]
    \Statex{\textbf{INPUT:}} Data points \(\{(x_i, y_i)\}\), number of cores \(c\), number of nodes \(M\), spline type (\texttt{cubic} or \texttt{b-spline}), spline degree \(k\)
    \Statex{\textbf{OUTPUT:}} $\T_{QTT}$, spline function \(s(x)\)
    \State \textbf{Sort} data points in increasing order based on \(x_i\)
    \State Set \(\text{start} \gets x_0\), \(\text{stop} \gets x_{-1}\)
    \If{\texttt{spline\_type} is \texttt{`cubic'}}
        \State \(s(x) \gets \text{CubicSpline}(\{(x_i, y_i)\})\)
    \ElsIf{\texttt{spline\_type} is \texttt{`b-spline'}}
        \State \(s(x) \gets \text{BSpline}(\{(x_i, y_i)\},\ \text{degree } k)\)
    \EndIf
    \State $\T_{QTT} \gets \texttt{interpolation\_qtt}(s(x),\ c,\ M)$
    \State \Return \(\T_{QTT},\ s(x)\)
  \end{algorithmic}
\end{algorithm}

By integrating the interpolated spline function $s(x)$ into the QTT framework, we effectively embed the empirical data into the tensor representation. This process enables the QTT model to learn from the data by capturing the essential features and patterns present in the observations. The resulting QTT tensor can then be utilized in solving PDEs, incorporating data-driven insights directly into the computational process. This approach enhances the model's accuracy and generalization capabilities, similar to how PINNs leverage data to inform their solutions. By learning from data, the QTT framework becomes more adaptable to complex real-world problems where data availability plays a crucial role.

\begin{figure}[ht]
    \centering
    \includegraphics[width=\linewidth]{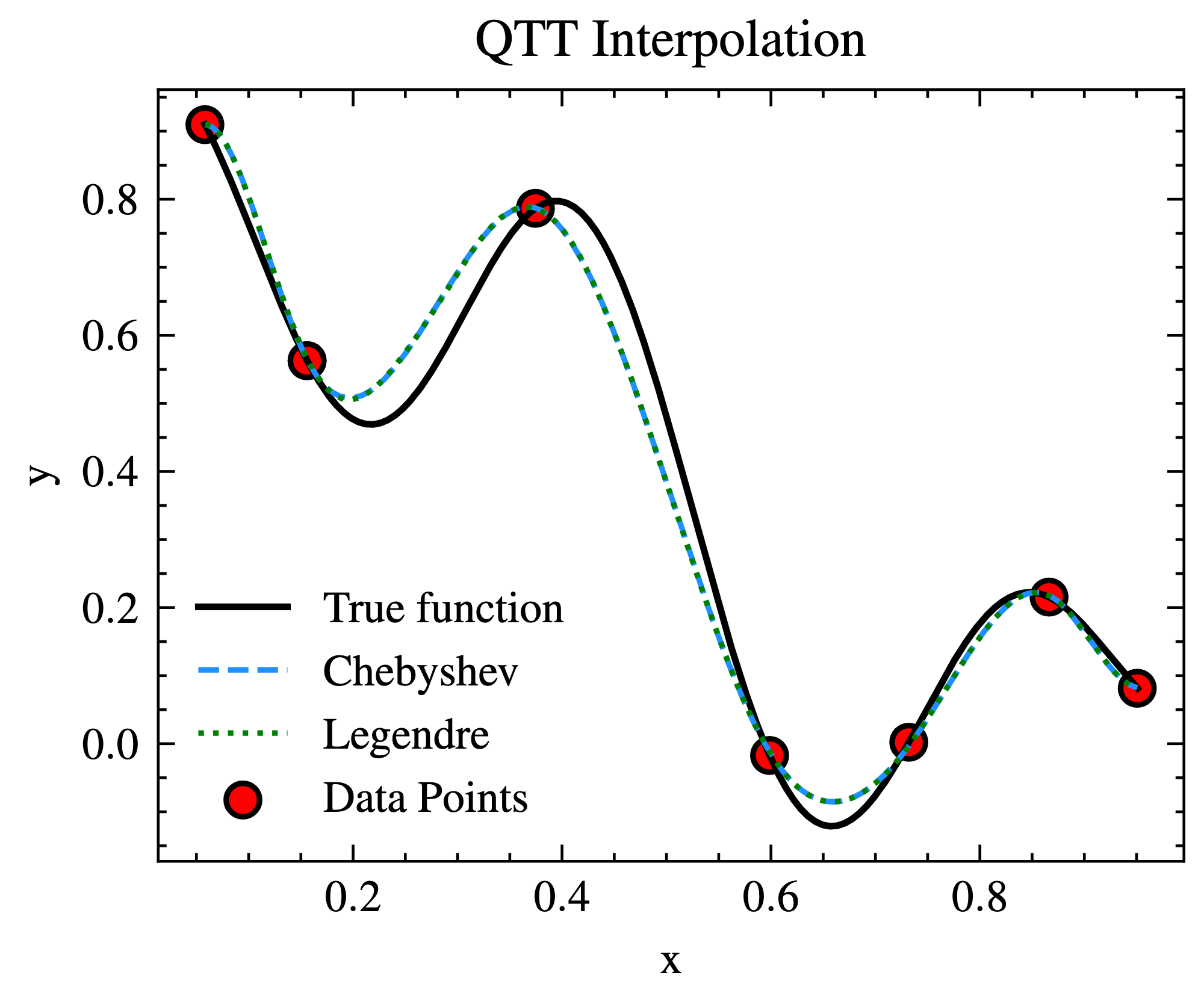}
    \caption{QTT interpolation on 7 data points sampled from a function $f(x)=\sin(3x)^2+\cos(5x)^3$. We used 8 cores and 25 interpolation nodes.} 
    \label{fig:one}
\end{figure}

An example of QTT interpolation on splines can be seen on figure \ref{fig:one} where the data points are sampled from the function $f(x)=\sin(3x)^2+\cos(5x)^3$ for both Chebyshev and Legendre nodes. From Figure \ref{fig:one} it is clear that a relatively low bond dimension is sufficient to capture the behavior of the function.

\subsection{Application to the Poisson Equation - Data-Driven Source Term}

To demonstrate the effectiveness of our data-driven approach within the QTT framework, we consider the following experiment. We begin with the Poisson equation:

\begin{align*}
    \nabla^2 u &= f(x,y) =  2x(y-1)(y-2x+xy+2)e^{x-y},
\end{align*}
where $(x,y) \in (0,1) \times (0,1)$, subject to all zero Dirichlet boundary conditions. Our goal is to learn a QTT representation of the source term $f$ using randomly sampled data points. We then solve the PDE using this learned QTT representation and compare the result to the analytical solution:

\begin{align*}
    u(x,y) = x(1-x)y(1-y)e^{x-y}.
\end{align*}

The first step of this experiment involves adapting Algorithm~\ref{alg:qtt_spline_interpolation} to use a bivariate spline and a 2D interpolation scheme. Once the QTT representation of the learned function is constructed, we employ the same algorithm from Section~\ref{sec:fd_poisson} to solve the Poisson equation, replacing the exact QTT representation of the source term with the learned QTT representation. Table \ref{tab:qtt_poisson_data} presents results for different discretization levels and numbers of data points. The column labeled ``Best $<$ 1 sec'' shows results obtained using MALS with an appropriate interpolation rank, while the ``MSE $<$ e-04'' column reports results from ALS tuned for higher speed. Compared to solving the same PDE without incorporating data points (see Table \ref{tab:problem_3}), we observe only a slight increase in runtime. This demonstrates that our method can maintain high accuracy while keeping the flexibility to trade off speed versus precision.

\begin{table}[h]
    \centering
    \begin{tabular*}{0.48\textwidth}{@{\extracolsep{\fill}}cccccccc@{}}
        \toprule
        \multicolumn{1}{c}{\textbf{Cores }} & \multicolumn{1}{c}{\textbf{\#Data}} & \multicolumn{2}{c}{\textbf{Best $<$ 1 sec}} & \multicolumn{2}{c}{\textbf{MSE $<$ e-04}}\\ 
        \cmidrule(lr){1-2} \cmidrule(lr){3-4} \cmidrule(lr){5-6}
        \textbf{p/ dim}    & \textbf{Points}  & \textbf{Time(s)} & \textbf{MSE} & \textbf{Time(s)} & \textbf{MSE}\\
        \hline
        \midrule
        10  & 64   & 0.589 & 7.66e-08 & 0.005 & 1.73e-04  \\
        10  & 128  & 0.803 & 3.52e-08 & 0.007 & 1.61e-04 \\
        10  & 256  & 0.764 & 3.86e-08 & 0.008 & 1.62e-04 \\
        12  & 256  & 0.712 & 1.02e-07 & 0.008 & 4.25e-04 \\
        14  & 256  & 0.655 & 8.41e-08 & 0.009 & 4.89e-04 \\
        \hline
        \bottomrule
    \end{tabular*}
    \caption{}
    \label{tab:qtt_poisson_data}
\end{table}

\begin{figure}[ht]
    \centering
    \includegraphics[width=\linewidth]{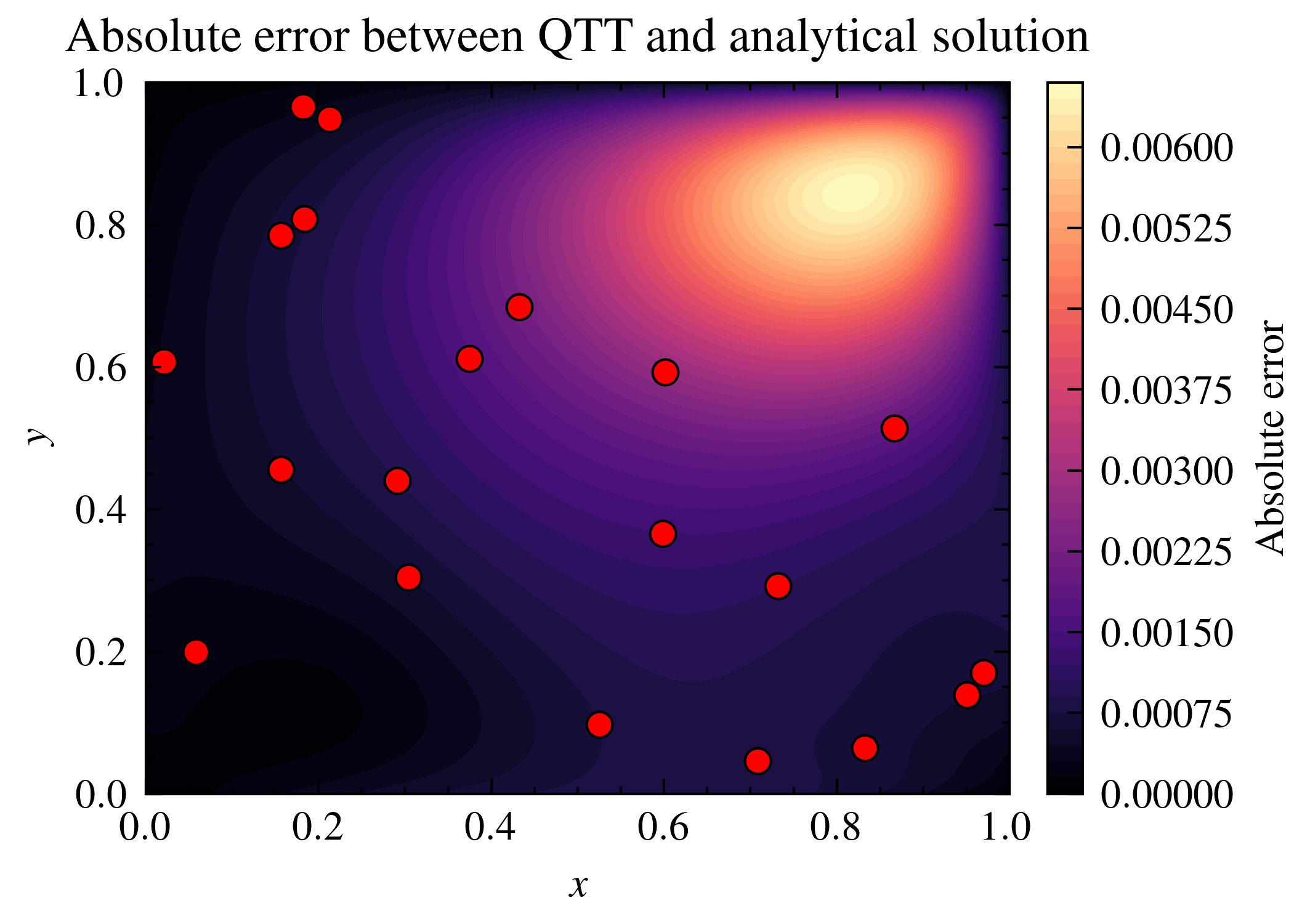}
    \caption{Point-wise absolute error between the QTT and the analytical solution, with an average MSE of order $10^{-6}$.}
    \label{fig:poisson_data_error}
\end{figure}

Figure \ref{fig:poisson_data_error}, shows the point-wise absolute error between the QTT solution with 10 cores per dimension and the analytical solution, where the source term was learned using 20 data points.

These results highlight that, in terms of both accuracy and speed, our QTT-based framework operates at a fundamentally different level than PINNs. For instance, Tables 8 and 12 from \cite{pinns_bench} show that for a similarly complex 2D Poisson equation the most accurate PINN configuration achieves an MSE of the order of $10^{-5}$, but requires approximately 800 seconds of runtime and their fastest solver is on the order of 300 seconds for an MSE of $10^{-4}$. This suggests that while PINNs might excel at rapid prototyping, the QTT framework dramatically outperforms ML based methods in terms of precision/speed.

The code is available on \url{https://github.com/DIKU-Quantum/TT-PDE}.

\section{Discussion}
We introduce a quantum-inspired solver that combines high speed with broad applicability for partial differential equations. Classical discretizations such as finite-element or multigrid methods, and even modern physics-informed neural networks, scale polynomially with the number of grid points. By contrast, our algorithm compresses both operators and solution fields into low-rank quantized tensor-train (QTT) formats. This compression reduces memory and runtime to logarithmic complexity relative to the underlying grid while preserving spectral-level accuracy. The same compressed framework also admits observational data as additional constraints, allowing data-driven refinement without a separate training phase.

Looking ahead, the most compelling next step is to test our framework on truly high-dimensional, real-world PDEs—settings in which conventional solvers become prohibitively expensive. Prime targets include (i) multi-asset option pricing in quantitative finance, (ii) high-dimensional Fokker–Planck and related kinetic equations in physics, and (iii) turbulence-resolved simulations in computational fluid dynamics. Success in these domains would showcase the solver’s ability to cut memory and runtime costs by orders of magnitude while retaining high fidelity.

Beyond these immediate applications, our space-time formulation also offers a fresh perspective on strongly nonlinear problems such as the Navier–Stokes equations. Future research will examine how low-rank QTT representations behave in chaotic regimes, develop a rigorous convergence theory that provides error bounds and stability guarantees for nonlinear operators, and quantify the bond dimensions needed to capture multiscale features without over-parameterising the model. Taken together, these studies will tighten the theoretical underpinnings of the method while extending its reach across scientific machine learning and computational science.

\paragraph{Acknowledgments}
We thank Nikita Gourianov, Sebastian Loeschke, and Alan Engsig-Karup for helpful discussions. 
We acknowledge support from the Carlsberg foundation and the Novo Nordisk foundation. 

\bibliography{main}

~ \newpage ~ \newpage 
\onecolumngrid

\appendix

\section{Useful QTT Constructions}
\label{ap:Usefull_QTT_Const}

In this section, we present key constructions for building QTT representation of matrices which are utilized across the solvers discussed in the main paper.

\subsection{Build QTT representation of ``Boundary Vector''}\label{appen:1d_bc_qtt}

Given the \textit{boundary vector}
\(\textbf{v} = \begin{pNiceMatrix}v_a & 0 & \cdots & 0 & v_b \end{pNiceMatrix}^\intercal\)
of length $2^c$, $c \geq 3$ the following construction builds the QTT representation of $\textbf{v}$ with $c$ cores:
\begin{align*}
\mathbf{v}=\F^{1,2,1,2}
 \bowtie
  \M1^{2,2,1,2} \bowtie
  (\M2^{(2,2,1,2)})^{\bowtie (c-3)} \bowtie
  \mathrm{L}^{(2,2,1,1)},
\end{align*}
with all the entries equal to zero except: ${\F_{0,0,0,1}=v_a},~  {\F_{0,0,0,1}=v_b},~{\M1_{0,1,0,0}=-1,}~{\M1_{1,0,0,1}=1},~{\M2_{0,1,0,0}=1},$ ${\M2_{1,0,0,1}=1}$ and ${\mathrm{L}_{1,0,0,0}=1},~{\mathrm{L}_{0,1,0,0}=-1}$

\subsection{Build ``Diagonal QTT'' from vector}\label{appen:vec_diag_mpo}

Given \(\textbf{v} = \begin{pNiceMatrix}v_1 & v_2 & \cdots & v_{2^c} \end{pNiceMatrix}^\intercal\) as input this method builds the $c$ cores MPO representation of ${D = \textrm{diag}(v_1,v_2,\hdots,v_{2^c})}$. The first step is:

\begin{center}
    \tikzset{every picture/.style={line width=0.75pt}}         

\begin{tikzpicture}[x=0.75pt,y=0.75pt,yscale=-1,xscale=1]

%Straight Lines [id:da6698841491839167] 
\draw [line width=0.75]    (193,21) -- (266.5,21) ;
\draw [shift={(268.5,21)}, rotate = 180] [color={rgb, 255:red, 0; green, 0; blue, 0 }  ][line width=0.75]    (10.93,-3.29) .. controls (6.95,-1.4) and (3.31,-0.3) .. (0,0) .. controls (3.31,0.3) and (6.95,1.4) .. (10.93,3.29)   ;
%Shape: Circle [id:dp614519652979976] 
\draw   (284.5,20.5) .. controls (284.5,13.32) and (290.32,7.5) .. (297.5,7.5) .. controls (304.68,7.5) and (310.5,13.32) .. (310.5,20.5) .. controls (310.5,27.68) and (304.68,33.5) .. (297.5,33.5) .. controls (290.32,33.5) and (284.5,27.68) .. (284.5,20.5) -- cycle ;
%Shape: Circle [id:dp12325261877394922] 
\draw   (332.5,20.5) .. controls (332.5,13.32) and (338.32,7.5) .. (345.5,7.5) .. controls (352.68,7.5) and (358.5,13.32) .. (358.5,20.5) .. controls (358.5,27.68) and (352.68,33.5) .. (345.5,33.5) .. controls (338.32,33.5) and (332.5,27.68) .. (332.5,20.5) -- cycle ;
%Shape: Circle [id:dp26194293166257143] 
\draw   (413.5,20.5) .. controls (413.5,13.32) and (419.32,7.5) .. (426.5,7.5) .. controls (433.68,7.5) and (439.5,13.32) .. (439.5,20.5) .. controls (439.5,27.68) and (433.68,33.5) .. (426.5,33.5) .. controls (419.32,33.5) and (413.5,27.68) .. (413.5,20.5) -- cycle ;
%Straight Lines [id:da3060029194908327] 
\draw    (297.5,33.5) -- (297.5,49.5) ;
%Straight Lines [id:da20349135052042133] 
\draw    (345.5,33.5) -- (345.5,49) ;
%Straight Lines [id:da019688301008758402] 
\draw    (426.5,33.5) -- (426.5,49) ;
%Straight Lines [id:da5811794437367376] 
\draw    (310.5,20.5) -- (332.5,20.5) ;
%Straight Lines [id:da10843409733863252] 
\draw    (358.5,20.5) -- (373.25,20.5) ;
%Straight Lines [id:da23415172023847552] 
\draw    (398.75,20.5) -- (413.5,20.5) ;
%Straight Lines [id:da932803402001576] 
\draw  [dash pattern={on 0.84pt off 2.51pt}]  (378.5,20.5) -- (393.25,20.5) ;

% Text Node
\draw (25,10.4) node [anchor=north west][inner sep=0.75pt]    {$\textbf{v}\ =\ ( v_{1} \ \ v_{2} \ \cdots \ v_{2^{c}})^{\intercal}$};
% Text Node
\draw (195,4) node [anchor=north west][inner sep=0.75pt]   [align=left] {MPS of $\displaystyle v$};
\draw (192,24) node [anchor=north west][inner sep=0.75pt]   [align=left] {Any Method};
\draw (220,38) node [anchor=north west][inner sep=0.75pt]   [align=left] {from};
\draw (197,52) node [anchor=north west][inner sep=0.75pt]   [align=left] {Section \ref{sec:qtt_rep_fun}};
% Text Node
\draw (290,14) node [anchor=north west][inner sep=0.75pt]    {$V_{1}$};
% Text Node
\draw (338,14) node [anchor=north west][inner sep=0.75pt]    {$V_{2}$};
% Text Node
\draw (420,14) node [anchor=north west][inner sep=0.75pt]    {$V_{c}$};

\end{tikzpicture}
\end{center}

\noindent Next, we build the following MPO:

\begin{center}
\tikzset{every picture/.style={line width=0.75pt}} %set default line width to 0.75pt        

\begin{tikzpicture}[x=0.75pt,y=0.75pt,yscale=-1,xscale=1]

%Shape: Circle [id:dp614519652979976] 
\draw   (24.9,42.5) .. controls (24.9,35.32) and (30.72,29.5) .. (37.9,29.5) .. controls (45.08,29.5) and (50.9,35.32) .. (50.9,42.5) .. controls (50.9,49.68) and (45.08,55.5) .. (37.9,55.5) .. controls (30.72,55.5) and (24.9,49.68) .. (24.9,42.5) -- cycle ;
%Shape: Circle [id:dp12325261877394922] 
\draw   (72.9,42.5) .. controls (72.9,35.32) and (78.72,29.5) .. (85.9,29.5) .. controls (93.08,29.5) and (98.9,35.32) .. (98.9,42.5) .. controls (98.9,49.68) and (93.08,55.5) .. (85.9,55.5) .. controls (78.72,55.5) and (72.9,49.68) .. (72.9,42.5) -- cycle ;
%Shape: Circle [id:dp26194293166257143] 
\draw   (153.9,42.5) .. controls (153.9,35.32) and (159.72,29.5) .. (166.9,29.5) .. controls (174.08,29.5) and (179.9,35.32) .. (179.9,42.5) .. controls (179.9,49.68) and (174.08,55.5) .. (166.9,55.5) .. controls (159.72,55.5) and (153.9,49.68) .. (153.9,42.5) -- cycle ;
%Straight Lines [id:da3060029194908327] 
\draw    (37.9,55.5) -- (37.9,71.5) ;
%Straight Lines [id:da20349135052042133] 
\draw    (85.9,14) -- (85.9,29.5) ;
%Straight Lines [id:da019688301008758402] 
\draw    (166.9,55.5) -- (166.9,71) ;
%Straight Lines [id:da5811794437367376] 
\draw    (50.9,42.5) -- (72.9,42.5) ;
%Straight Lines [id:da10843409733863252] 
\draw    (98.9,42.5) -- (113.65,42.5) ;
%Straight Lines [id:da23415172023847552] 
\draw    (139.15,42.5) -- (153.9,42.5) ;
%Straight Lines [id:da932803402001576] 
\draw  [dash pattern={on 0.84pt off 2.51pt}]  (118.9,42.5) -- (133.65,42.5) ;
%Straight Lines [id:da9576701627045947] 
\draw    (37.9,13.5) -- (37.9,29.5) ;
%Straight Lines [id:da34897053807319567] 
\draw    (85.9,55.5) -- (85.9,71.5) ;
%Straight Lines [id:da8443240574800431] 
\draw    (166.9,13.5) -- (166.9,29.5) ;
%Shape: Circle [id:dp5072902725266553] 
\draw   (253.3,42.9) .. controls (253.3,35.72) and (259.12,29.9) .. (266.3,29.9) .. controls (273.48,29.9) and (279.3,35.72) .. (279.3,42.9) .. controls (279.3,50.08) and (273.48,55.9) .. (266.3,55.9) .. controls (259.12,55.9) and (253.3,50.08) .. (253.3,42.9) -- cycle ;
%Straight Lines [id:da7911925459164109] 
\draw    (266.3,14.4) -- (266.3,29.9) ;
%Straight Lines [id:da5236643265452031] 
\draw    (231.3,42.9) -- (253.3,42.9) ;
%Straight Lines [id:da5955944269511714] 
\draw    (279.3,42.9) -- (294.05,42.9) ;
%Straight Lines [id:da8125342173264719] 
\draw    (266.3,55.9) -- (266.3,71.9) ;
%Shape: Circle [id:dp08640849141175477] 
\draw   (358.5,29.7) .. controls (358.5,22.52) and (364.32,16.7) .. (371.5,16.7) .. controls (378.68,16.7) and (384.5,22.52) .. (384.5,29.7) .. controls (384.5,36.88) and (378.68,42.7) .. (371.5,42.7) .. controls (364.32,42.7) and (358.5,36.88) .. (358.5,29.7) -- cycle ;
%Straight Lines [id:da7382247077929036] 
\draw    (336.5,29.7) -- (358.5,29.7) ;
%Straight Lines [id:da34534753670241547] 
\draw    (384.5,29.7) -- (399.25,29.7) ;
%Straight Lines [id:da2969152508741423] 
\draw    (371.5,42.7) -- (371.5,58.7) ;

% Text Node
\draw (29,35) node [anchor=north west][inner sep=0.75pt]    {$D_{1}$};
% Text Node
\draw (77,35) node [anchor=north west][inner sep=0.75pt]    {$D_{2}$};
% Text Node
\draw (158,35) node [anchor=north west][inner sep=0.75pt]    {$D_{c}$};
% Text Node
\draw (190.4,35) node [anchor=north west][inner sep=0.75pt]   [align=left] {s.t.};
% Text Node
\draw (414,23) node [anchor=north west][inner sep=0.75pt]   [align=left]{
$\text{if } i=j,~ \forall \alpha ,\beta,$
};
% Text Node
\draw (314.4,4) node [anchor=north west][inner sep=0.75pt]    {$\begin{cases}
 & \\
 & \\
 & \\
 & 
\end{cases}$};
% Text Node
\draw (237.6,30) node [anchor=north west][inner sep=0.75pt]    {$\alpha $};
% Text Node
\draw (255.8,35) node [anchor=north west][inner sep=0.75pt]    {$D_{k}$};
% Text Node
\draw (281.6,27) node [anchor=north west][inner sep=0.75pt]    {$\beta $};
% Text Node
\draw (268,7.8) node [anchor=north west][inner sep=0.75pt]    {$j$};
% Text Node
\draw (268,56.6) node [anchor=north west][inner sep=0.75pt]    {$i$};
% Text Node
\draw (299,38.9) node [anchor=north west][inner sep=0.75pt]    {$=$};
% Text Node
\draw (342.8,17) node [anchor=north west][inner sep=0.75pt]    {$\alpha $};
% Text Node
\draw (363,22) node [anchor=north west][inner sep=0.75pt]    {$V_{k}$};
% Text Node
\draw (386.8,14) node [anchor=north west][inner sep=0.75pt]    {$\beta $};
% Text Node
\draw (373.2,43.4) node [anchor=north west][inner sep=0.75pt]    {$i$};
% Text Node
\draw (365,65.4) node [anchor=north west][inner sep=0.75pt]    {$0$};
% Text Node
\draw (414,64.2) node [anchor=north west][inner sep=0.75pt]   [align=left] {otherwise.};
\end{tikzpicture}
\end{center}

\subsection{Build ``Eraser QTT''}\label{appen:eraser_qtt}
Given the constants $n_1,n_2,n_3,n_4$ and $c$ this method construct a QTT representation of $\mathrm{M}$ a $2^c \times 2^c$ matrix with all its elements equal to zero except for the four entries shown below:
\begin{align*}
\mathrm{M} = \begin{bNiceMatrix} 
n_1 & n_2 & \phantom{0} & \phantom{0} & \phantom{\dots}  \\ 
\phantom{0} & 0 & 0 & \phantom{0} & \phantom{\dots}  \\  
&  &  \Ddots & \Ddots   \\ 
& &  &  0  & 0   \\
& &  &  & n_3 & n_4 
\end{bNiceMatrix}_{(2^c \times 2^c)},
\end{align*}
the construction is given by:
\begin{align*}
\mathrm{M} = \F^{1,2,2,2}
 \bowtie
  (\M^{(2,2,2,2)})^{\bowtie (c-2)} \bowtie
  \mathrm{L}^{(2,2,1,1)},
\end{align*}
with all the entries equal to zero except for: $\F_{0,0,0,0}=1$, $\F_{0,1,1,1}=1$,  
$\M_{0,0,0,0}=1$, $\M_{1,1,1,1}=1$, $\mathrm{L}_{0,0,0,0}=n_1$, ${\mathrm{L}_{0,0,1,0}=n_2}$, $\mathrm{L}_{1,1,0,0}=n_3$ and $\mathrm{L}_{1,1,1,0}=n_4$ 

\subsection{Building the Analytic QTT Representation of $f(x) = \sin(\alpha x + \phi)$}\label{appen:analytic_qtt_sine}

We construct an analytic rank-2 QTT representation of $f(x) = \sin(\alpha x + \phi)$, where $x$ is discretized in the interval $(0,1)$ with $2^c$ grid points. The trick of this construction is to use the trigonometric identity ${\sin(\alpha \pm \beta) = \sin(\alpha)\cos(\beta) \pm \cos(\alpha) \sin(\beta)}$. Using our indexing convention, the QTT representation of the discretized function $f$ is given by:
\begin{align*}
\F^{1,2,1,2}
 \bowtie
  (\M^{(2,2,1,2)})^{\bowtie (c-2)} \bowtie
  \mathrm{L}^{(2,2,1,1)},
\end{align*}
where the tensor components are defined as follows:
\begin{align*}
    F^{0,0,0,:} &= (\sin(\phi), \cos(\phi)), \\
    F^{0,1,0,:} &= (\sin(\alpha x[2^{c-1}] + \phi), \cos(\alpha x[2^{c-1}] + \phi)),
\end{align*}
for the middle cores:
\begin{align*}
    M^{0,:,0,0} &= (1, \cos(\alpha x[2^i]))^\intercal, 
    \quad M^{0,:,0,1} = (0, -\sin(\alpha x[2^i]))^\intercal, \\
    M^{1,:,0,0} &= (0, \sin(\alpha x[2^i]))^\intercal,
    \quad M^{1,:,0,1} = (1, \cos(\alpha x[2^i]))^\intercal,
\end{align*}
for $i = c-2,\ldots,0$ and the last core is given by:
\begin{align*}
    L^{0,:,0,0} &= (1, \cos(\alpha x[1]))^\intercal, \\
    L^{1,:,0,0} &= (0, \sin(\alpha x[1]))^\intercal.
\end{align*}

\subsubsection{Adjusting to a shifted interval}

While the above construction assumes $x$ is discretized in $(0,1)$ with $2^c$ points, our solvers typically require values on the shifted discrete interval $(h, 1-h)$, where $h = \frac{1}{2^c + 1}$. To transition to this target grid, we proceed as follows: First, discretize $(0,1)$ into $2^c + 2$ points and define the target sequence: $y_i = \frac{i}{2^c + 1}, \text{ for } i=1,2,\dots,2^c.$ Next, discretize $(0,1)$ into $2^c$ points and define the sequence: $x_i = \frac{i-1}{2^c - 1}, \text{ for } i=1,2,\dots,2^c.$ To map between these two grids, we use the transformation:
\begin{align*}
    \sin(\alpha x + \phi) = \sin(Ky),
\end{align*}
where $K$ is constant and holds for: $\alpha = K \frac{2^c - 1}{2^c + 1}, \text{ and } \phi = \frac{K}{2^c + 1}.$ This transformation ensures that whenever we need the QTT representation of $f(x)$ on the adjusted interval, we can obtain it via an appropriate rescaling of the parameters.\\

\subsubsection{Higher Dimensions}

Extending the analytical QTT representation from the 1D to the 2D function ${f(x,y) = \sin(\alpha_1 x + \phi_1)\sin(\alpha_2 y + \phi_2)}$
is straightforward since the 2D function can be represented as the tensor (Kronecker) product of each 1D component. In practice, we first construct the QTT representation of each component with the desired number of cores and interval and then concatenate these 1D components to get the QTT representation of the full discretized 2D function in serial ordering. This same construction naturally extends to higher-dimensions.

\subsection{Building the Analytic QTT Representation of $f(x) = e^{\alpha x}$}\label{appen:analytic_qtt_exp}

We construct the analytic rank-1 QTT representation of $f(x) = e^{\alpha x}$, where $x$ is discretized in the interval $(0,1)$ with $2^c$ grid points. Using our indexing convention, the QTT representation of the discretized function $f$ function is given by:
\begin{align*}
  \F_1^{(1,2,1,1)} \bowtie \F_2^{(1,2,1,1)} \bowtie \cdots \bowtie \F_c^{(1,2,1,1)}
\end{align*}
where the tensor components are defined as:
\begin{align*}
    F_i^{0,:,0,0} &= 
    \begin{bmatrix} 
    1 \\ \exp(\alpha x[2^{c-i}])
    \end{bmatrix}.
\end{align*}
We note that the same procedure described in the previous section can be applied to extend this construction to higher-dimensions.

\section{QTT Solvers for PDEs}
\label{ap:Other_Solvers}
\subsection*{1D Heat Equation Time Stepping vs Space-time}\label{appen:heat_ts_st}
In this section, we first construct a standard time-stepping finite difference scheme for the 1D heat equation in the QTT framework. We then develop a QTT space-time formulation of the same problem. To evaluate their efficiency, we compare the runtime and accuracy of both approaches against the analytical solution of the heat equation for a specific test case.

We consider the one-dimensional heat equation:
\begin{align}
      &\frac{\partial u}{\partial t} =  \frac{\partial^2 u}{\partial x^2},  \quad (x,t) \in \Omega = \{a \leq x \leq b, t \geq 0\},\label{heat_1_pde}
\end{align}
with initial condition $u(x,0) = g_0(x)$ and boundary conditions $u(a,t)=g_1(t)$ and $u(b,t)=g_2(t)$.

Let $w_{ij} \approx u(x_i,t_j)$ and consider the following discretization of \eqref{heat_1_pde}:
\begin{align}\label{heat_1_disc}
    \frac{1}{l}(w_{i,j+1} - w_{i,j}) =  
    \frac{1}{h^2}(w_{i-1,j+1} - 2w_{i,j+1} + w_{i+1,j+1}),
\end{align}
where $l$ is the time step size and $h = (b-a)/(N+1)$, with $N=2^c$, the dimension time step. We can express \eqref{heat_1_disc} as a system of linear equations $A w_{j+1} = w_j + b_{j+1}$, with
\begin{align*}
A =\begin{bNiceMatrix}
1-2r & r &  &  &    \\
r & 1-2r & r &  &    \\
& \ddots & \ddots & \ddots     \\
\end{bNiceMatrix}_{(N \times N)}, 
w_j = 
\begin{bNiceMatrix}
    w_{1,j+1} \\ w_{2,j+1}\\ \Vdots \\ w_{N,j+1}
\end{bNiceMatrix}_{(N \times 1)},
b_{j+1} = 
\begin{bNiceMatrix}
    r w_{0,j+1} \\ 0 \\ \Vdots \\ 0 \\ r w_{N,j+1}
\end{bNiceMatrix}_{(N \times 1)},
\end{align*}
with $r = -l/h^2$. The main steps of the algorithm to solve this differential equation are given below:

\begin{algorithm}[H]
\caption{QTT Solver 1D Heat Equation}\label{algorithm_1d_heat}
\begin{algorithmic}[1]
\Statex{INPUT:} $c$, $l$, \textit{timesteps}

\Statex{OUTPUT:} $c$ cores QTT representation of the solution at each time step

\Statex{Build $c$ cores QTT representation of:}

\State $\textrm{A}_{QTT}$ \Comment{Lemma \eqref{lemma_1}}

\State $w_0$ from initial conditions \Comment{Any method from \ref{sec:qtt_rep_fun}}

\For{$k=0$ until (\textit{timesteps}-1)}

\State Build QTT rep of $b_{k+1}$ \Comment{Using Construction \ref{appen:1d_bc_qtt}}
\State $w_{sol}[k+1] \leftarrow$ALS($\textrm{A}_{QTT}, w_{sol}[k],w_{sol}[k]+(4)$)

\EndFor
\State Return{} $w_{sol}$
\end{algorithmic}
\end{algorithm}
\noindent This algorithm is highly efficient since all steps in the main loop are performed entirely within the tensor framework. The primary computational cost comes from running ALS over the required number of time steps.

Now, we apply a space-time discretization by treating the time variable $t$ as an additional spatial dimension. Consider the one-dimensional heat equation given by \eqref{heat_1_pde} but now $t$ is defined over the same region as $x$.
Let $w_{ij} \approx u(t_i,x_j)$ and consider the following discretization of \eqref{heat_1_pde}:
\begin{align}\label{heat1d_space}
    \frac{1}{h}(w_{i-1,j} - w_{i,j}) - \frac{1}{h^2} (-w_{i,j-1} + 2w_{i,j} - w_{i,j+1}) = 0,
\end{align}
where $h = (b-a)/N$, with $N=2^c$, representing the discretization step. We can express \eqref{heat1d_space} as the system of linear equations $Aw=b$, where
\begin{align*}
A=\frac{1}{h}\begin{pmatrix}
-1 & 0 &  &  &    \\
1 & -1 & 0 &  &    \\
 & 1 & -1 & 0  &    \\
  & & \ddots & \ddots & \ddots     \\
\end{pmatrix}_{2^c \times 2^c} 
\otimes 
\I_{2^c} 
-
\I_{2^c}
\otimes
\frac{1}{h^2}\begin{pmatrix}
2 & -1 &  &  &    \\
-1 & 2 & -1 &  &    \\
 & -1 & 2 & -1  &    \\
 & & \ddots & \ddots & \ddots    \\
\end{pmatrix}_{2^c \times 2^c},
\end{align*}
\begin{align}\label{w_spacetime}
    w &= \begin{bmatrix} w_{1,1} & w_{1,2} &\cdots & w_{1,N} & w_{2,1}& \cdots & w_{2,N} & \cdots\cdots  & w_{N,N} \end{bmatrix}^T,
\end{align}
and the right-hand side $b$ encodes the initial and boundary conditions:
\begin{align*}
    b = -1 (| 0 \rangle \otimes (w_{0,x})/h + (w_{t,0})/h^2 \otimes | 0 \rangle + (w_{t,1})/h^2 \otimes | 1 \rangle),
\end{align*}
the $| 0 \rangle$ and $| 1 \rangle$ are the vectors of length $N$ of the form $(1,0,\hdots,0)$ and $(0,0,\hdots,1)$, respectively.
The space-time algorithm is very similar to the one presented in Section \ref{sec:fd_poisson} to solve Poisson equation. By Lemma \eqref{lemma_1}, we can construct a low-rank QTT representation of $A$, while the QTT representation of $b$ can be computed using any of the discussed methods in Section \ref{sec:qtt_rep_fun}. Once these components are built, we apply (M)ALS to obtain a QTT representation of the solution. This solution can be interpreted as performing $2^c$ time steps of size $1/2^c$ from the previous algorithm.

To compare both methods, we consider the one-dimensional heat equation on the domain $0 \leq x \leq 1$ with the initial condition: $$u(x,0) = \sin \left( \frac{\pi x}{2} \right) + \frac{1}{2} \sin (2 \pi x)$$ and boundary conditions: $$u(0,t)=0, u(1,t)=\exp(-\pi^2 t/4).$$ The analytical solution is given by: $$u(x,t) = \exp(- \pi^2 t / 4) \sin \left(\frac{\pi x}{2}\right) + \frac{1}{2} \exp(-4\pi^2 t) \sin (2\pi x).$$

Table \ref{tab:spacetime_heat} reports the results for the space-time QTT solver using MALS with a single sweep.
\begin{table}[h]
    \centering
    \begin{tabular*}{0.4\textwidth}{@{\extracolsep{\fill}}ccc@{}}
        \toprule
        \textbf{Cores p/ dim} & \textbf{Run Time (s)} & \textbf{MSE} \\ 
        \hline
        \midrule
        6  & 0.00463  & 8.27e-05  \\  
        8  & 0.00515  & 6.10e-06  \\  
        10  & 0.00686  & 3.99e-07    \\  
        12  & 0.00812  & 2.51e-08    \\  
        14  & 0.01041  & 2.2e-09     \\  
        \hline
        \bottomrule
    \end{tabular*}
    \caption{}
    \label{tab:spacetime_heat}
\end{table}

For comparison, Table \ref{tab:timestepping} presents the results for the time-stepping QTT solver optimized to obtain the same MSE.
\begin{table}[h]
    \centering
    \begin{tabular*}{0.5\textwidth}{@{\extracolsep{\fill}}cccc@{}}
        \toprule
        \textbf{Cores p/ dim} & \textbf{Time Steps} & \textbf{Run Time (s)} & \textbf{Avg MSE} \\ 
        \hline
        \midrule
        4   & $2^5$   & 0.04245  & 9.68273e-05  \\  
        5  & $2^7$   & 0.20569  & 8.3352e-06   \\  
        6  & $2^{9}$ & 0.93338  & 5.706e-07    \\  
        7  & $2^{10}$ & 2.21900  & 7.17e-08     \\  
        8  & $2^{11}$ & 4.99420  & 8.9e-09      \\  
        \hline
        \bottomrule
    \end{tabular*}
    \caption{}
    \label{tab:timestepping}
\end{table}

The time step size is given by $1/\#\text{(time steps)}$. The results demonstrate that the space-time method is significantly more efficient than the traditional time-stepping approach. As the number of cores increases, the space-time method improves accuracy with a minimal runtime growth. In contrast, the time step method is more accurate regarding the number of discretization points, but still requires a significant amount of cores and time steps to achieve high-accuracy solutions.

\subsection*{2D Heat Equation with time-dependent boundary conditions}\label{appen:heat_2d_td_bc}
In this section, we analyze how our QTT framework handles the 2D heat equation with time-dependent boundary conditions, focusing on algorithmic aspects and the efficient treatment of these conditions. A key advantage of combining QTT with interpolation techniques is the ability to incorporate complex, time-varying boundaries at each time step without significantly increasing the runtime. At the end of this section, we present a table summarizing the performance of the solver for a different number of discretization points and time steps.

We consider the 2D heat equation:
\begin{align}\label{heat_2d_pde}
\frac{\partial u}{\partial t} = \alpha \left( \frac{\partial^2 u}{\partial x^2} + \frac{\partial^2 u}{\partial y^2} \right), \end{align}
where $(x,y) \in \Omega = (0,1) \times (0,1)$, $t \geq 0$ and the thermal diffusivity $\alpha=0.6$.

For the left and top boundaries, we impose time-dependent boundary conditions modeled as two Gaussian sources (double Gaussian waveforms) moving along the edges:
\begin{align*} 
u(0,y,t) &= \frac{1}{\sqrt{2 \pi}} \left( e^{-10 (y + 2 - t)^2} + e^{-10 (y - 3.4 + t)^2} \right),\\
u(x,1,t) &= \frac{1.5}{\sqrt{2 \pi}} \left( e^{-10 (x + 2 - t)^2} + e^{-10 (x - 3.4 + t)^2 } \right). 
\end{align*}
The remaining boundary conditions and the initial condition are set to zero.

Physically, each boundary condition can be interpreted as two moving heat sources tracing linear trajectories but in opposite directions along the borders of a square. This setup allows us to test the solver's ability to handle dynamic boundary conditions efficiently while maintaining accurate results.

Let $w^k_{ij} \approx u(x_i,y_j,t_k)$ and consider the following discretization of \eqref{heat_2d_pde}:

\begin{align*}
    \frac{1}{l}(w^{k+1}_{i,j} - w^k_{i,j}) = 
    &\frac{1}{h^2}(w^k_{i-1,j+1} - 2w^k_{i,j+1} + w^k_{i+1,j+1}\\
    &+ w^k_{i+1,j-1} - 2w^k_{i+1,j} + w^k_{i+1,j+1}),
\end{align*}
where $l$ is the time step size and $h = (b-a)/(N+1)$, with $N=2^c$, the dimension time step. We can express the equation above as a system of linear equations $A w^{k+1} = w^k + b^{k+1}$, with $A = (H^{(c)} \otimes \I_{2^c} + \I_{2^c} \otimes H^{(c)})$, where
\begin{align}\label{H_m}
H^{(c)} =\begin{bNiceMatrix}
1-2r & r &  &  &    \\
r & 1-2r & r &  &    \\
& \ddots & \ddots & \ddots     \\
\end{bNiceMatrix}_{(N \times N)}, 
\end{align}
with $r = -l/h^2$ and the vectors $w^k$ and $b^{k+1}$ similarly to Equations \eqref{w_2d} and \eqref{b_vec} respectively.

The main steps of the algorithm to solve this PDE are given below:

\begin{algorithm}[H]
\caption{QTT Solver Heat Equation}\label{alg:heat_td_bc}
\begin{algorithmic}[1]
\Statex{INPUT:} $c$, $l$, \textit{timesteps}

\Statex{OUTPUT:} $c$ cores QTT representation of the solution at each time step

\Statex{Build $c$ cores QTT representation of:}

\State $\textrm{A}_{QTT}$ \Comment{Lemma \eqref{lemma_1}}

\State $w^0$ from initial conditions \Comment{Interpolation}

\For{$k=0$ until (\textit{timesteps}-1)}

\State Build QTT rep of $b^{k+1}$ \Comment{Interpolation}
\State $w_{sol}[k+1] \leftarrow$ALS($\textrm{A}_{QTT}, w_{sol}[k],w_{sol}[k]+(4)$)

\EndFor
\State Return{} $w_{sol}$
\end{algorithmic}
\end{algorithm}

To build the QTT representation of $b^{k+1}$, at each run of line $(4)$, we use a 1D rank-revealing interpolation scheme (Section \ref{sec:interpolation}) since one of the spatial dimensions will always be fixed and $t$ will be given by the constant $k \cdot l$. The time step size is always fixed as $1/\#\text{timesteps}$.

\begin{table}[H]
    \centering
    \begin{tabular*}{0.8\textwidth}{@{\extracolsep{\fill}}cccc@{}}
        \toprule
        \textbf{Cores p/ dim} & \textbf{Time Steps} & \textbf{Time(s) Build B.C.} & \textbf{Total Run Time(s)} \\ 
        \hline
        \midrule
        5   & 100    & 0.0717  & 0.30160  \\  
        5   & 1000   & 0.7127  & 2.97216  \\  
        5   & 10000  & 7.1249  & 29.69824 \\  
        10  & 100    & 0.1631  & 0.71565  \\  
        10  & 1000   & 1.6534  & 7.20218  \\  
        10  & 10000  & 16.9066 & 72.74729 \\
        \hline
        \bottomrule
    \end{tabular*}
    \caption{}
    \label{tab:V}
\end{table}

In Table \ref{tab:V} the column ``Time(s) Build B.C.''  corresponds to line (4) of the algorithm, executed for the total number of time steps.  We observe that computing the time-dependent boundary conditions with an interpolation scheme accounts for less than $25\%$ of the total runtime. Additionally, solving the PDE on a $2^{10} \times 2^{10}$ grid takes only about 2.5 times longer compared to a $2^{5} \times 2^{5}$ grid that is 1024 times smaller, highlighting the efficiency of our approach.

The table below shows the time required to construct the classical representation of $b^{k+1}$ (line (4) of the algorithm)
\begin{table}[H]
    \centering
    \begin{tabular*}{0.6\textwidth}{@{\extracolsep{\fill}}ccc@{}}
        \toprule
        \textbf{Cores p/ dim} & \textbf{Time Steps} & \textbf{Time(s) Build B.C.} \\
        \hline
        \midrule
        10  & 100    & 1.5393   \\  
        10  & 1000   & 15.4698  \\  
        10  & 10000  & 152.0644 \\  
        \hline
        \bottomrule
    \end{tabular*}
    \caption{}
    \label{}
\end{table}
Even if a classical solver could match the performance of the QTT solver, the time required to construct the time-dependent boundary conditions would still be the dominant computational cost, making the classical approach less efficient than our method.

\begin{figure}[H]
    \centering
    \includegraphics[width=0.99\linewidth]{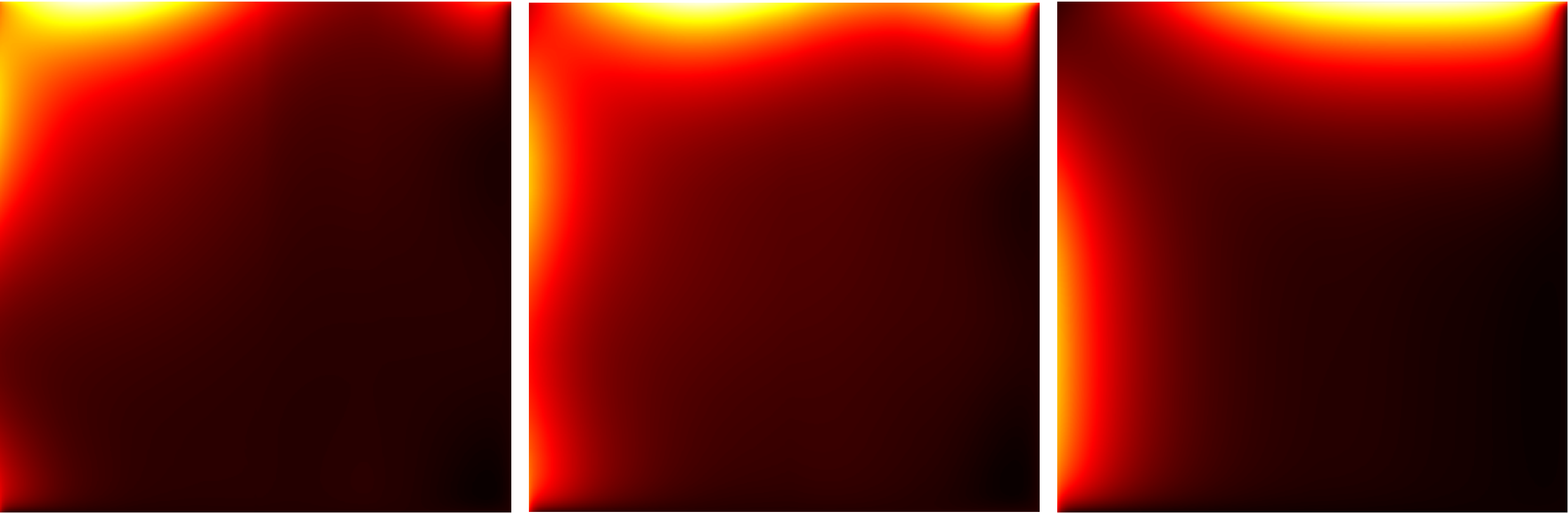}
    \caption{Solution of the 2D heat equation at different time steps.}
\end{figure}

~\newpage
\subsection*{1D Burgers' Equation Space-time Discretization Scheme}\label{appen:burgers_st}
We present the discretization scheme used to build Algorithm \ref{alg:burgers_spacetime} in the main paper to solve the one-dimensional Burgers' equation: 
\begin{align*}
      &\frac{\partial u}{\partial t} =  \nu \frac{\partial^2 u}{\partial x^2} - u\frac{\partial u}{\partial x}, 
\end{align*}
using a space-time approach. Let
$(t,x) \in [a,b] \times [c,d]$, with the initial condition $u(0,x) = g_0(x)$, and boundary conditions: $u(t,a)=g_1(t)$ and $u(t,b)=g_2(t)$.

For simplicity, assume $a=c=0$ and $b=d=1$. Let $w_{ij} \approx u(t_i,x_j)$ and consider the following discretization:
\begin{align}\label{burgers_spacetime_disc}
    \frac{w_{i-1,j} - w_{i,j}}{h} &- \nu
    \frac{- w_{i,j-1} + 2w_{i,j} - w_{i,j+1}}{h^2}
    +u\frac{w_{i,j-1} - w_{i,j+1}}{h}
    = 0,
\end{align}
where wlog $h = (b-a)/N$, with $N=2^c$, is the dimension step size. We can express \eqref{burgers_spacetime_disc} as the system of linear equations $Aw=b$, where
\begin{align*}
    A = \frac{1}{h}\begin{pmatrix}
-1 & 0 &  &  &    \\
1 & -1 & 0 &  &    \\
 & 1 & -1 & 0  &    \\
  & & \ddots & \ddots & \ddots     \\
\end{pmatrix}_{2^c \times 2^c} 
\otimes 
\I_{2^c} 
-
\I_{2^c}
\otimes
\frac{\nu}{h^2}\begin{pmatrix}
2 & -1 &  &  &    \\
-1 & 2 & -1 &  &    \\
 & -1 & 2 & -1  &    \\
 & & \ddots & \ddots & \ddots    \\
\end{pmatrix}_{2^c \times 2^c} 
    + u~  
\I_{2^c}
\otimes 
\frac{1}{h}
\begin{pmatrix}
0 & -1 &  &  &    \\
1 & 0 & -1 &  &    \\
& 1 & 0 & -1  &    \\
& & \ddots & \ddots & \ddots     \\
\end{pmatrix}_{2^c \times 2^c},
\end{align*}
and the right-hand side $b$ encodes the initial and boundary conditions:
\begin{align*}
    b = -1 (| 0 \rangle \otimes (w_{0,x})/h + (w_{t,0})/h^2 \otimes | 0 \rangle + (w_{t,1})/h^2 \otimes | 1 \rangle),
\end{align*}
the $| 0 \rangle$ and $| 1 \rangle$ are the vectors of length $N$ of the form $(1,0,\hdots,0)$ and $(0,0,\hdots,1)$, respectively. The vector $w$ follows the same indexing as the vector defined in Equation \eqref{w_spacetime}. The initial approximation of $u$ is given by a $2^c \times 2^c$ diagonal matrix $D_1$ whose elements correspond to the $2^c$ entries of the vector $b$. For the next \textit{runs}, the previous solution is used to refine the approximation of $u$.

\subsection*{Plot of the Solution of Burgers' Equation from Section \ref{sec:time_stepping_vs_ST}}\label{appen:burgers_plots}
Below, we present plots of the solution of Burgers' equation with parameters specified in Section~\ref{sec:time_stepping_vs_ST}:

\begin{center}
\begin{figure}[h]
    % First row (QTT solutions)
    \begin{minipage}{0.4\textwidth}
        \centering
        \includegraphics[width=\linewidth]{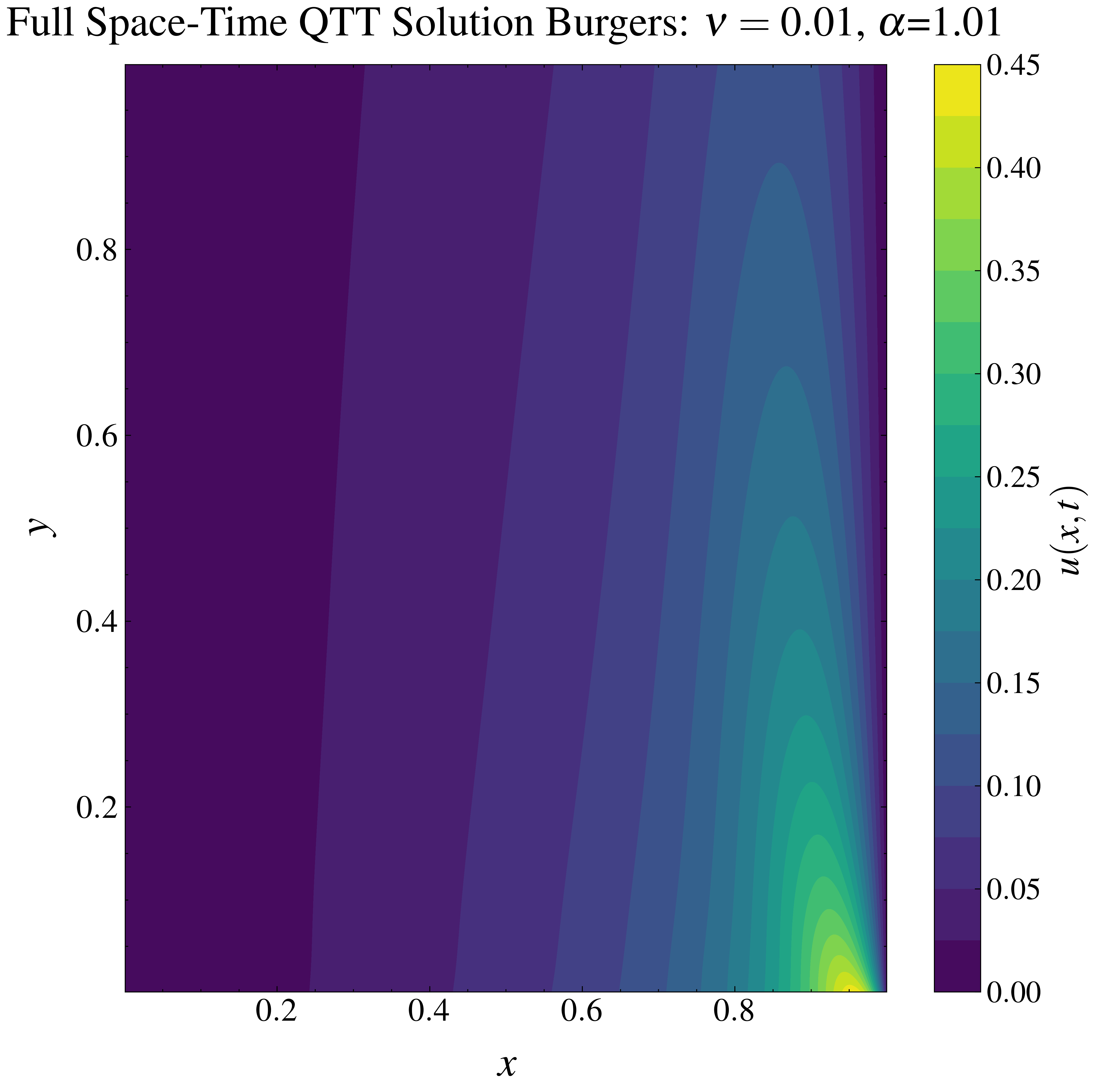}
    \end{minipage}
    %\hfill
    \begin{minipage}{0.4\textwidth}
        \centering
        \includegraphics[width=\linewidth]{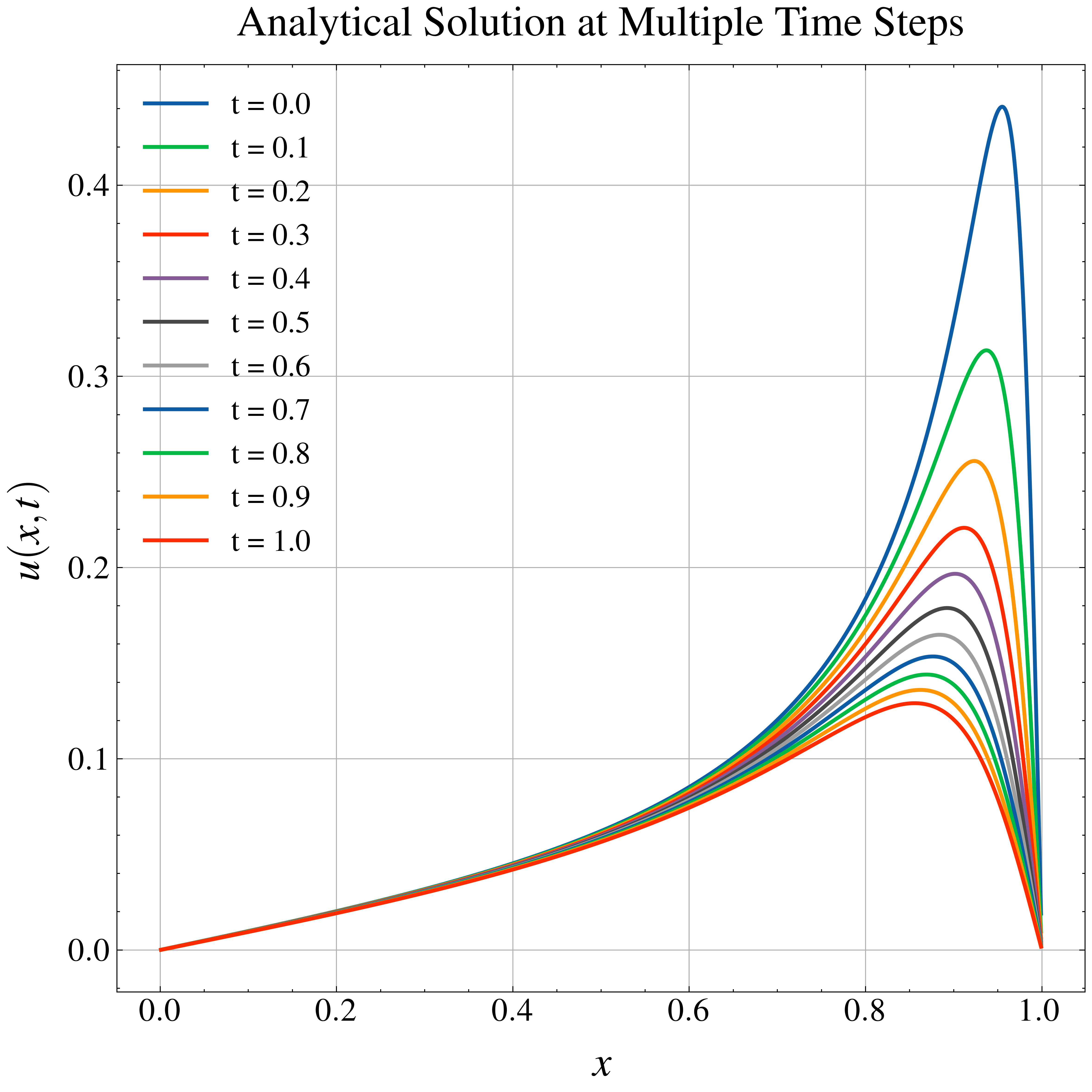}
    \end{minipage}
    \label{fig:burgers_solution}
\end{figure}
\end{center}

~\newpage
\section{Benchmarks}
\label{app:Benchmarks}
All experiments were conducted on a standard MacBook Pro with an Apple M3 Pro processor and 18GB of RAM. For problems 1 to 3 PyAMG, was configured in the following way:
\begin{verbatim}
    A = pyamg.gallery.poisson((2 ** c, 2 ** c), format='csr')  # Discrete Laplace 
    ml = pyamg.ruge_stuben_solver(A)          # construct the multigrid hierarchy
    sol = ml.solve(v_B, tol=1e-18)
\end{verbatim}
For the QTT solvers we used MALS with 2 sweeps and the initial guess of the solution was a random low-rank QTT. This seems to be a good configuration for a wide range of problems.

\subsection{Problem 1 - Section \ref{sec:Laplace_Poisson}}\label{appen:problem_1}
Here we have the results that were used to build the plot in Section \ref{sec:Laplace_Poisson} and also present a plot of the analytical solution. Consider the Laplace equation
\begin{align*}
    \nabla^2 u &= 0, \quad (x,y) \in (0,1) \times (0,1),
\end{align*}
with boundary conditions:
\begin{alignat*}{2}
    &u(x,0) = \sin(k \pi x) \sinh(k \pi), \quad &&u(x,1) = 0, \\
    &u(0,y) = 0, \quad  &&u(1,y) = 0.
\end{alignat*} The analytical solution is given by:
\begin{align*}
    u(x,y) = \sin(k \pi x) \sinh(k\pi(1-y)).
\end{align*}
The results in the last columns of Table \ref{tab_viii} were obtained with the QTT solver described in Section \ref{sec:Laplace_Poisson} where the boundary condition was built using an analytic QTT representation of the sine function (Section \ref{appen:analytic_qtt_exp}). 
\twocolumngrid
\noindent For $k=3$:
\begin{flushleft}
    \begin{table}[h]
    \begin{tabular*}{0.45\textwidth}{@{\extracolsep{\fill}}ccccc@{}}
        \toprule
        Cores & \multicolumn{2}{c}{\textbf{PyAMG}} & \multicolumn{2}{c}{\textbf{Analytic}} \\ 
        \cmidrule(lr){2-3} \cmidrule(lr){4-5}
        p/ dim & Time(s)   & MSE        & Time(s)   & MSE        \\ 
        \hline
        \midrule
        6  & 0.04166   & 1.51e+00   & 0.01281   & 1.61e+00   \\
        7  & 0.14485   & 9.91e-02   & 0.03262   & 1.02e-01   \\
        8  & 0.52393   & 6.34e-03   & 0.05292   & 6.44e-03   \\
        9  & 2.02764   & 4.01e-04   & 0.07461   & 4.04e-04   \\
        10 & 7.92881   & 2.52e-05   & 0.12971   & 2.53e-05   \\
        11 & 32.24839  & 1.58e-06   & 0.21816   & 1.58e-06   \\
        12 & 131.7988  & 9.89e-08   & 0.28115   & 9.91e-08   \\
        13 & -         & -          & 0.35427   & 6.19e-09   \\
        14 & -         & -          & 0.33853   & 3.82e-10   \\
        \hline
        \bottomrule
    \end{tabular*}
    \caption{}
    \label{tab_viii}
\end{table}
\end{flushleft}

\onecolumngrid

\begin{picture}(0,0)
    \put(300,0){\includegraphics[width=7cm]{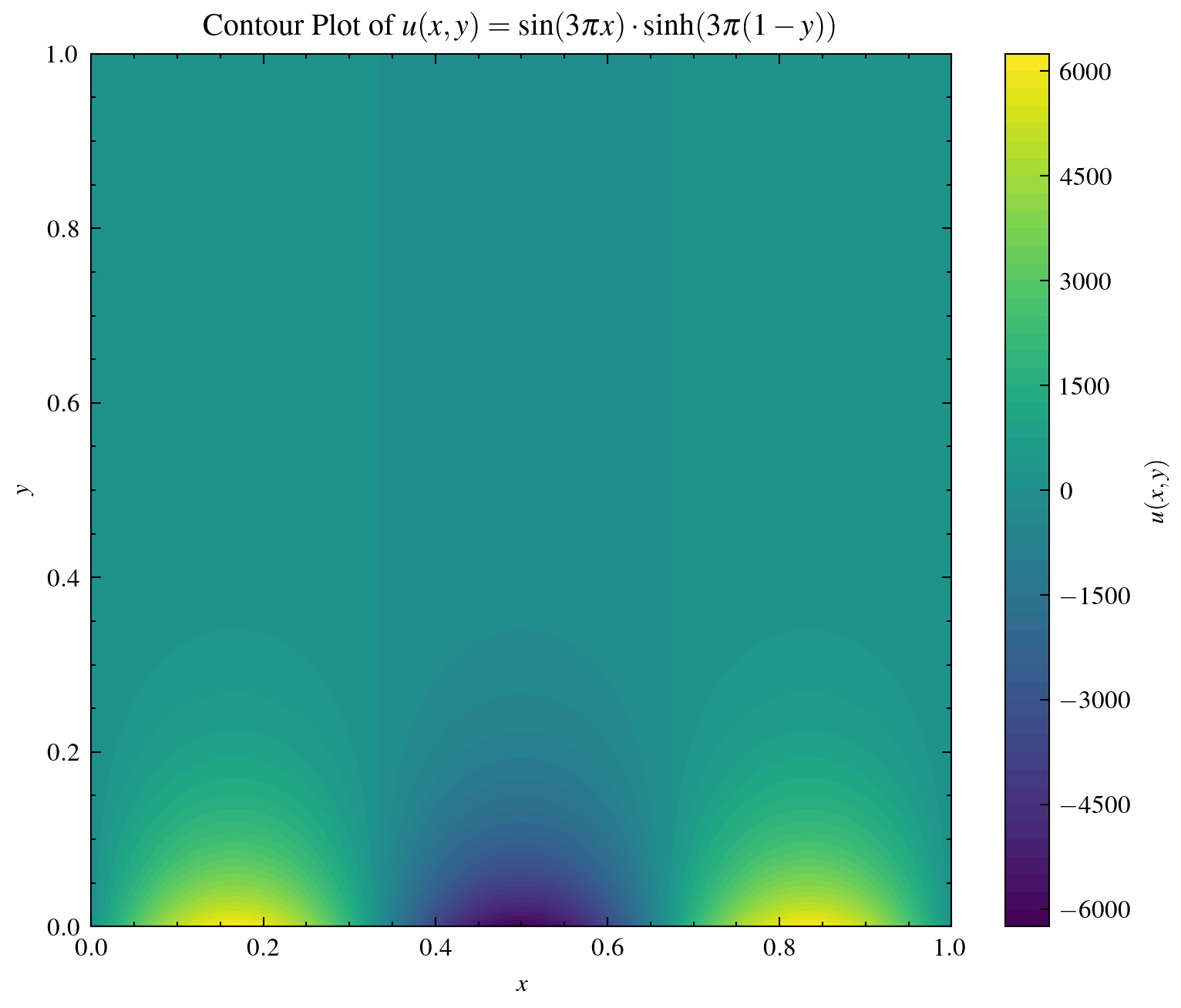}}
\end{picture}

\subsection{Problem 2 - Section \ref{sec:Laplace_Poisson}}\label{appen:problem2}
In Section \ref{sec:Laplace_Poisson} we consider the Poisson equation:
\begin{align*}
    \frac{\partial^2u}{\partial x^2} + \epsilon_1 \frac{\partial^2u}{\partial y^2} + \epsilon_2 \frac{\partial^2u}{\partial z^2}&= -\sin(\pi x) \sin(\pi y) \sin(\pi z), \quad (x,y,z) \in (0,1)^3,
\end{align*}
with all zero boundary conditions 
for $\epsilon_1=0.001$ and $\epsilon_2=0.0001$ now we present the results for $\epsilon_1$ and $\epsilon_2=1$. The analytical solution is given by:
\begin{align*}
    u(x,y) = \dfrac{\sin(\pi x) \sin(\pi y) \sin(\pi z)}{\pi^2 (1 +\epsilon_1 + \epsilon_2)}.
\end{align*}

\begin{table}[h]
    \begin{tabular*}{0.7\textwidth}{@{\extracolsep{\fill}}ccccccc@{}}
        \toprule
        Cores & \multicolumn{2}{c}{\textbf{PyAMG}} & \multicolumn{2}{c}{\textbf{Analytic}} & \multicolumn{2}{c}{\textbf{Interpolation}} \\ 
        p/ dim & Run Time(s)   & MSE        & Run Time(s)   & MSE        & Run Time(s)   & MSE\\ \hline
        2  & 0.01146   & 3.14e-07   & 0.00101   & 3.14e-07   & 0.00227   & 1.43e-05   \\
        3  & 0.01558   & 2.12e-08   & 0.00119   & 2.12e-08   & 0.01619   & 5.65e-06   \\
        4  & 0.05651   & 1.39e-09   & 0.00179   & 1.39e-09   & 0.08303   & 1.73e-06   \\
        5  & 0.42155   & 8.93e-11   & 0.00442   & 8.93e-11   & 0.26738   & 4.76e-07   \\
        6  & 3.15405   & 5.66e-12   & 0.00902   & 5.66e-12   & 0.42960   & 1.25e-07   \\
        7  & 24.93626  & 3.57e-13   & 0.01986   & 3.57e-13   & 0.42929   & 3.19e-08   \\
        8  & 215.11537 & 2.24e-14   & 0.02987   & 2.24e-14   & 0.19902   & 8.05e-09   \\
        9  & -         & -          & 0.03497   & 1.40e-15   & 0.13490   & 2.02e-09   \\
        \hline
    \end{tabular*}
    \caption{}
    \label{tab_ix}
\end{table}

Analyzing Table \ref{tab_ix} we note that the PyAMG and ``QTT Analytic solver'' successfully achieve the same grid discretization error, up to a rounding factor. The small variation in runtime among the QTT solvers is due to the randomness of the initial guess for the solution. The ``QTT Interpolation solver'' constructs the boundary condition using a rank-revealing method (see Section \ref{sec:qtt_rep_fun}), which requires more cores for higher accuracy but eventually become competitive with PyAMG.
\subsection{Problem 3}\label{appen:problem_3}
For the next problem we will consider other two types of QTT solver. The first solver uses the same default configuration but builds the source term using TT-SVD. The “QTT Optimized” solver has the same configuration as “QTT Interpolation,” but its initial random guess has smaller ranks than the former.
Consider Poisson equation:
\begin{align*}
    \nabla^2 u &= 2x(y-1)(y-2x+xy+2)e^{x-y}, \quad (x,y) \in (0,1) \times (0,1),
\end{align*}
with all zero boundary conditions. The analytical solution is given by:
\begin{align*}
    u(x,y) = x(1-x)y(1-y)e^{x-y}.
\end{align*}

\begin{table}[h]
    \begin{tabular*}{0.95\textwidth}{@{\extracolsep{\fill}}ccccccccc@{}}
        \toprule
        Cores & \multicolumn{2}{c}{\textbf{PyAMG}} & \multicolumn{2}{c}{\textbf{TT-SVD}} & \multicolumn{2}{c}{\textbf{Interpolation}} & \multicolumn{2}{c}{\textbf{Optimized}} \\ 
        p/ dim & Run Time(s)   & MSE        & Run Time(s)   & MSE        & Run Time(s)   & MSE        & Run Time(s)   & MSE\\ \hline
        5  & 0.02714   & 2.54e-10   & 0.01664   & 2.87e-10   & 0.04307   & 1.21e-07   & 0.02545   & 1.21e-07   \\
        6  & 0.04353   & 1.74e-11   & 0.02067   & 1.85e-11   & 0.09251   & 3.35e-08   & 0.04450   & 3.35e-08   \\
        7  & 0.14212   & 1.14e-12   & 0.03075   & 1.17e-12   & 0.14752   & 8.82e-09   & 0.05660   & 8.82e-09   \\
        8  & 0.50735   & 7.27e-14   & 0.05338   & 7.39e-14   & 0.20932   & 2.27e-09   & 0.07040   & 2.27e-09   \\
        9  & 1.96091   & 4.60e-15   & 0.08975   & 4.63e-15   & 0.22444   & 5.75e-10   & 0.06070   & 5.75e-10   \\
        10 & 8.01757   & 2.89e-16   & 0.14419   & 2.90e-16   & 0.15338   & 1.45e-10   & 0.07000   & 1.45e-10   \\
        11 & 32.36982  & 1.81e-17   & 0.31486   & 1.81e-17   & 0.22861   & 3.63e-11   & 0.07240   & 3.62e-11   \\
        12 & -         & -          & 0.85459   & 1.13e-18   & 0.19887   & 9.08e-12   & 0.07620   & 9.08e-12   \\
        13 & -         & -          & 3.11220   & 6.77e-20   & 0.19800   & 2.26e-12   & 0.06880   & 3.08e-12   \\
        \hline
    \end{tabular*}
    \caption{}
    \label{tab:problem_3}
\end{table}

We note that PyAMG and the ``QTT TT-SVD solver'' successfully achieve the grid discretization error and a simple modification on the ``QTT Interpolation solver'' manages to get the same accuracy but with a further 10 times speed-up.

\subsection{Problem 4}\label{appen:burgers_pinns}
\enlargethispage{0.3\baselineskip}
In this section, we analyze the effect of the \textit{runs} parameter in the Space-Time QTT Algorithm~\ref{alg:burgers_spacetime} on the solution of a specific instance of Burgers' equation. We consider the equation with the following initial and boundary conditions:
\begin{align*}\label{burgers_eq}
      &\frac{\partial u}{\partial t} =  (0.01/\pi) \frac{\partial^2 u}{\partial x^2} - u\frac{\partial u}{\partial x}, \quad x \in [-1,1], t \in [0,1]\\ 
    &u(x,0) = -\sin(\pi x),~ u(-1,t)=u(1,t)=0.
\end{align*}

The following plots were generated by running MALS with 2 sweeps and 10 cores in each spatial dimension starting from a random initial guess of the solution:
\begin{center}
\begin{figure}[h]
    % First row (QTT solutions)
    \begin{minipage}{0.24\textwidth}
        \centering
        \includegraphics[width=\linewidth]{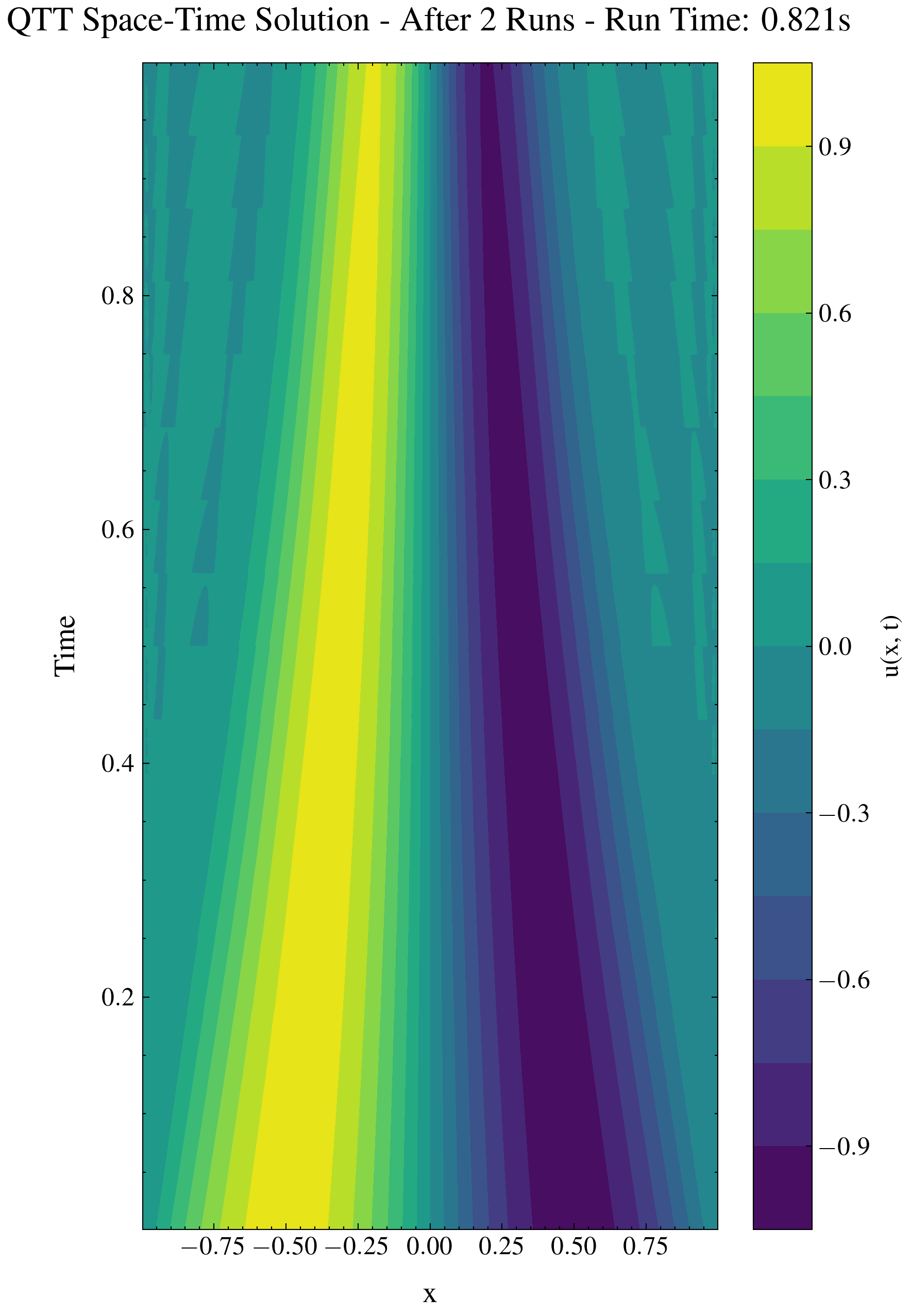}
    \end{minipage}
    \hfill
    \begin{minipage}{0.24\textwidth}
        \centering
        \includegraphics[width=\linewidth]{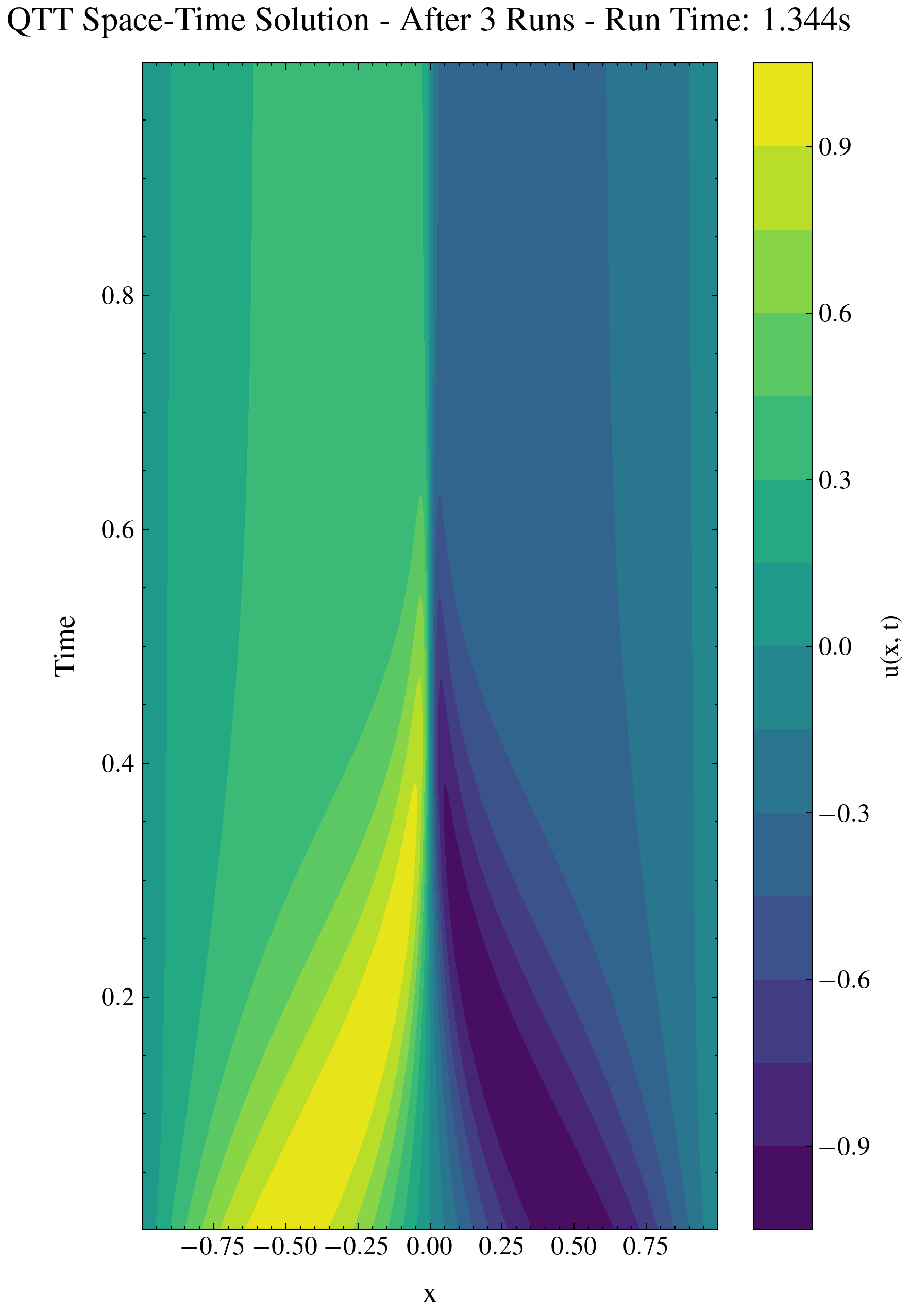}
    \end{minipage}
    \hfill
    \begin{minipage}{0.24\textwidth}
        \centering
        \includegraphics[width=\linewidth]{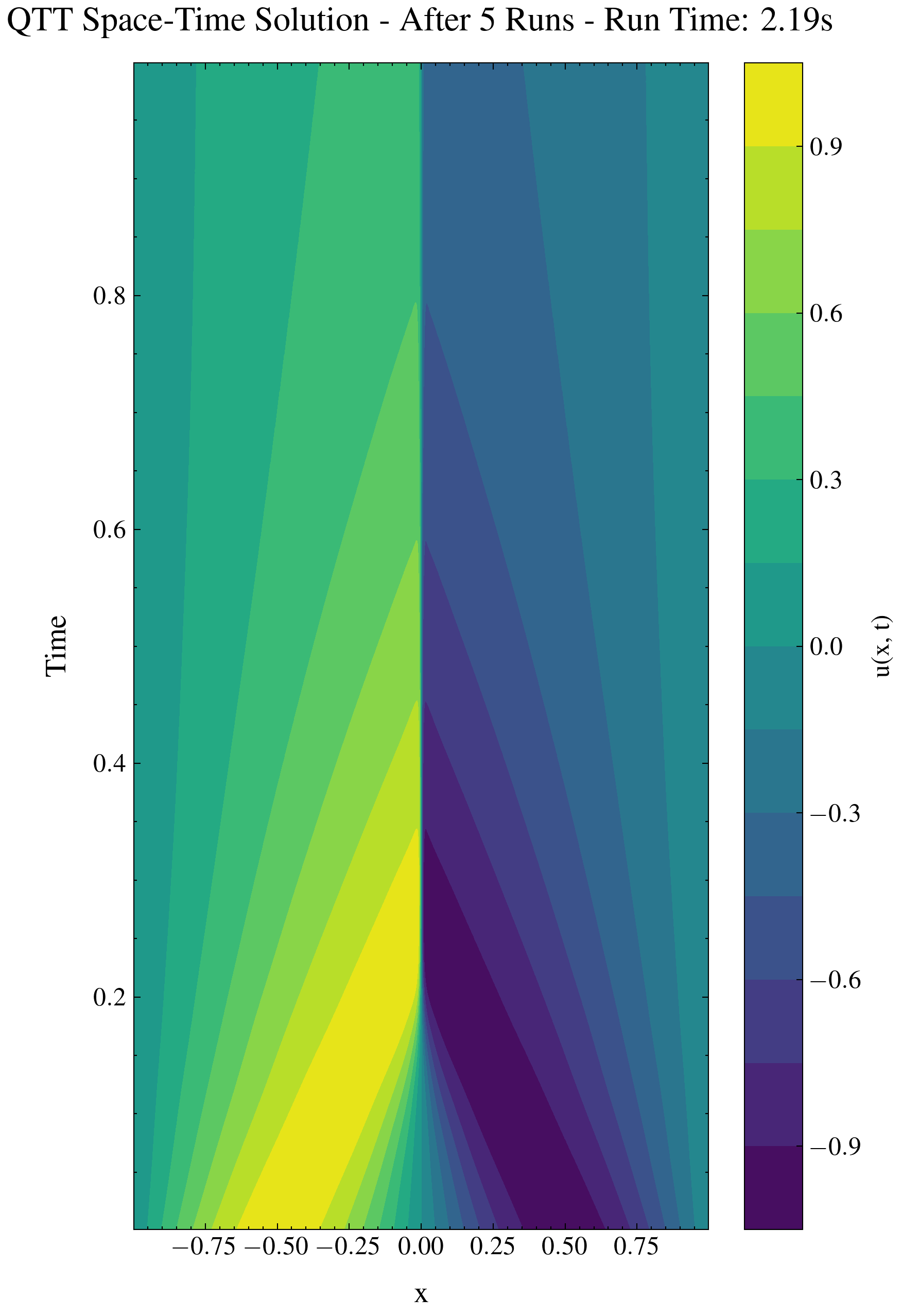}
    \end{minipage}
    \hfill
    \begin{minipage}{0.24\textwidth}
        \centering
        \includegraphics[width=\linewidth]{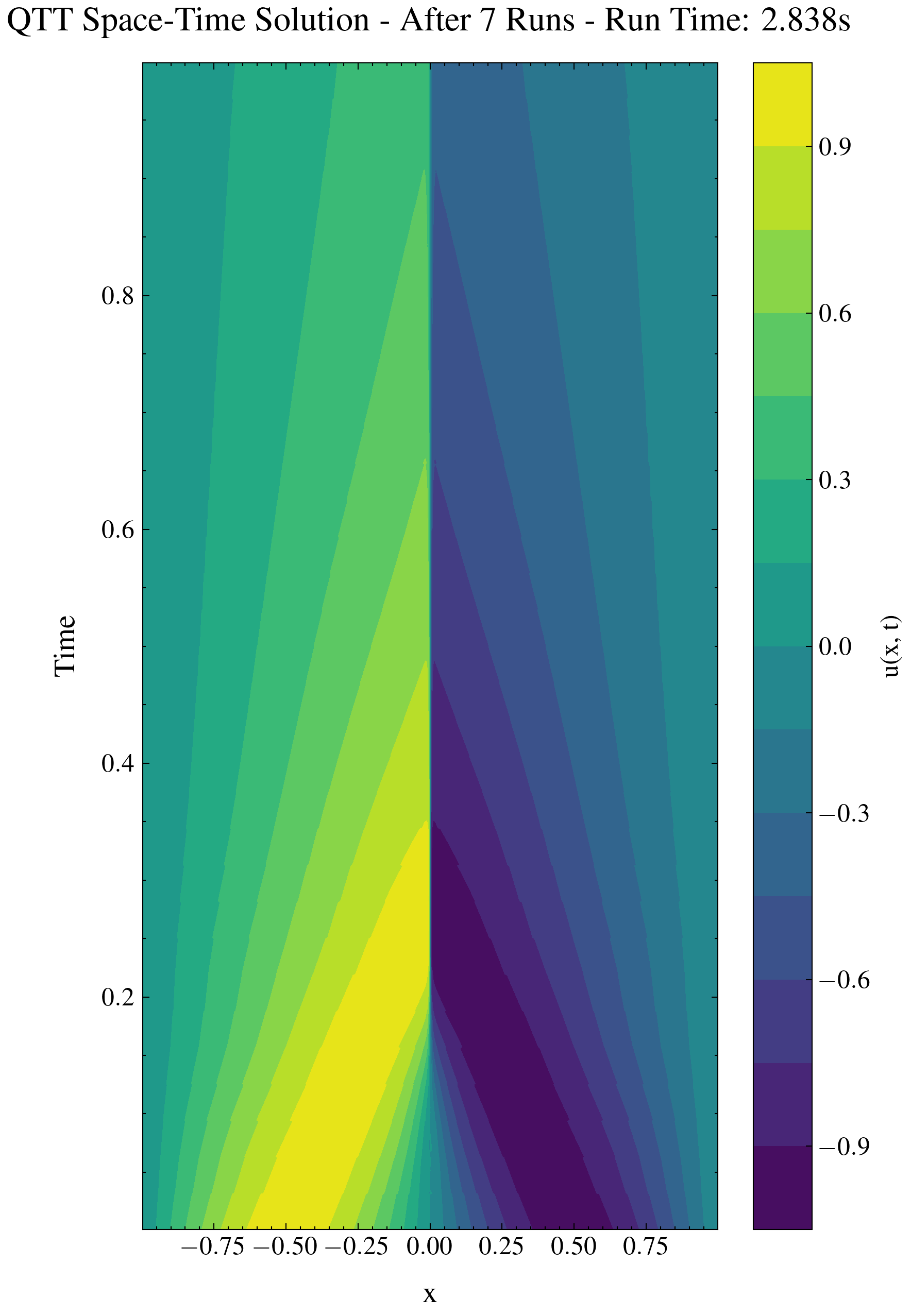}
    \end{minipage}

    \vspace{0.3cm} % Reduce spacing between rows

    % Second row (Slices)
    \begin{minipage}{0.24\textwidth}
        \centering
        \includegraphics[width=\linewidth]{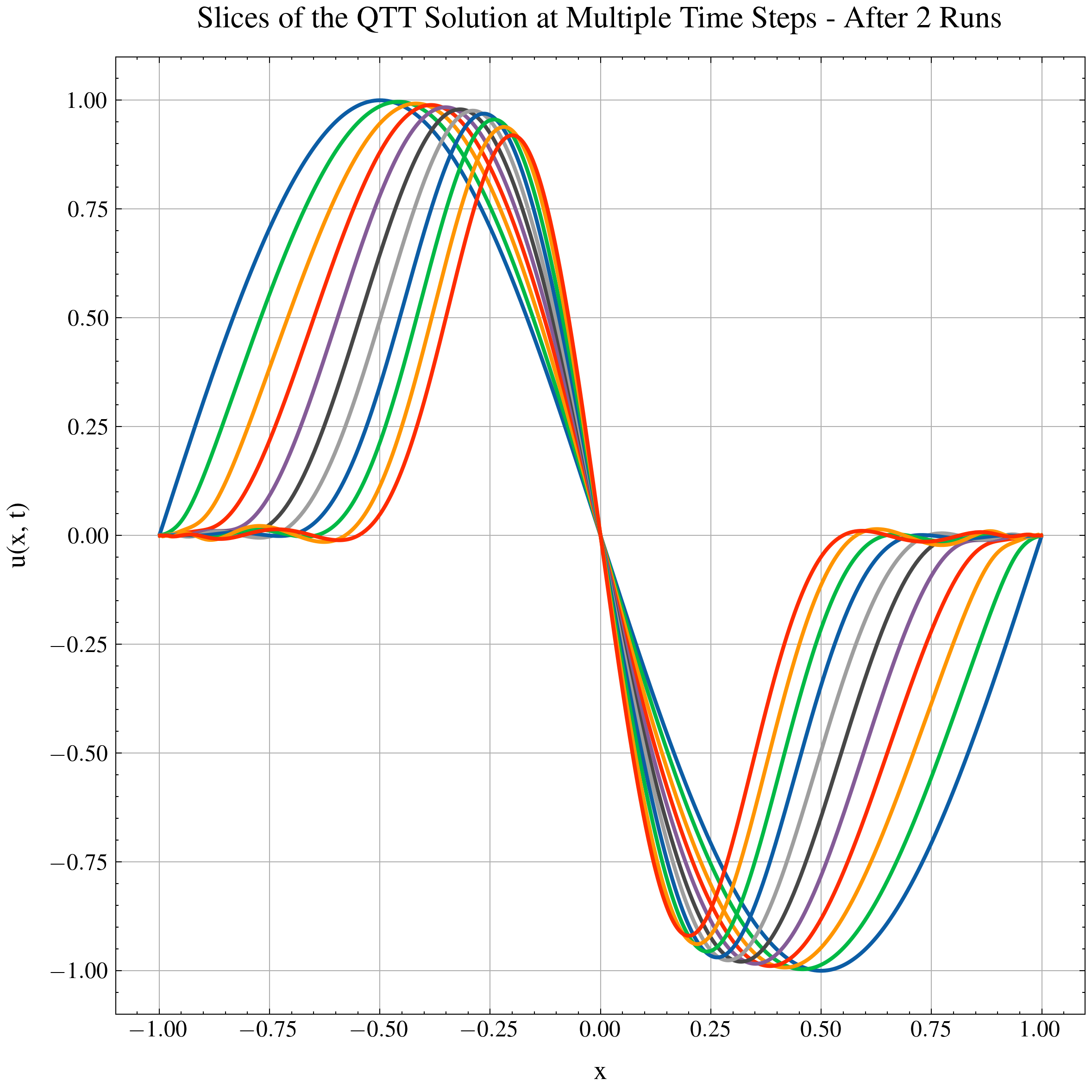}
    \end{minipage}
    \hfill
    \begin{minipage}{0.24\textwidth}
        \centering
        \includegraphics[width=\linewidth]{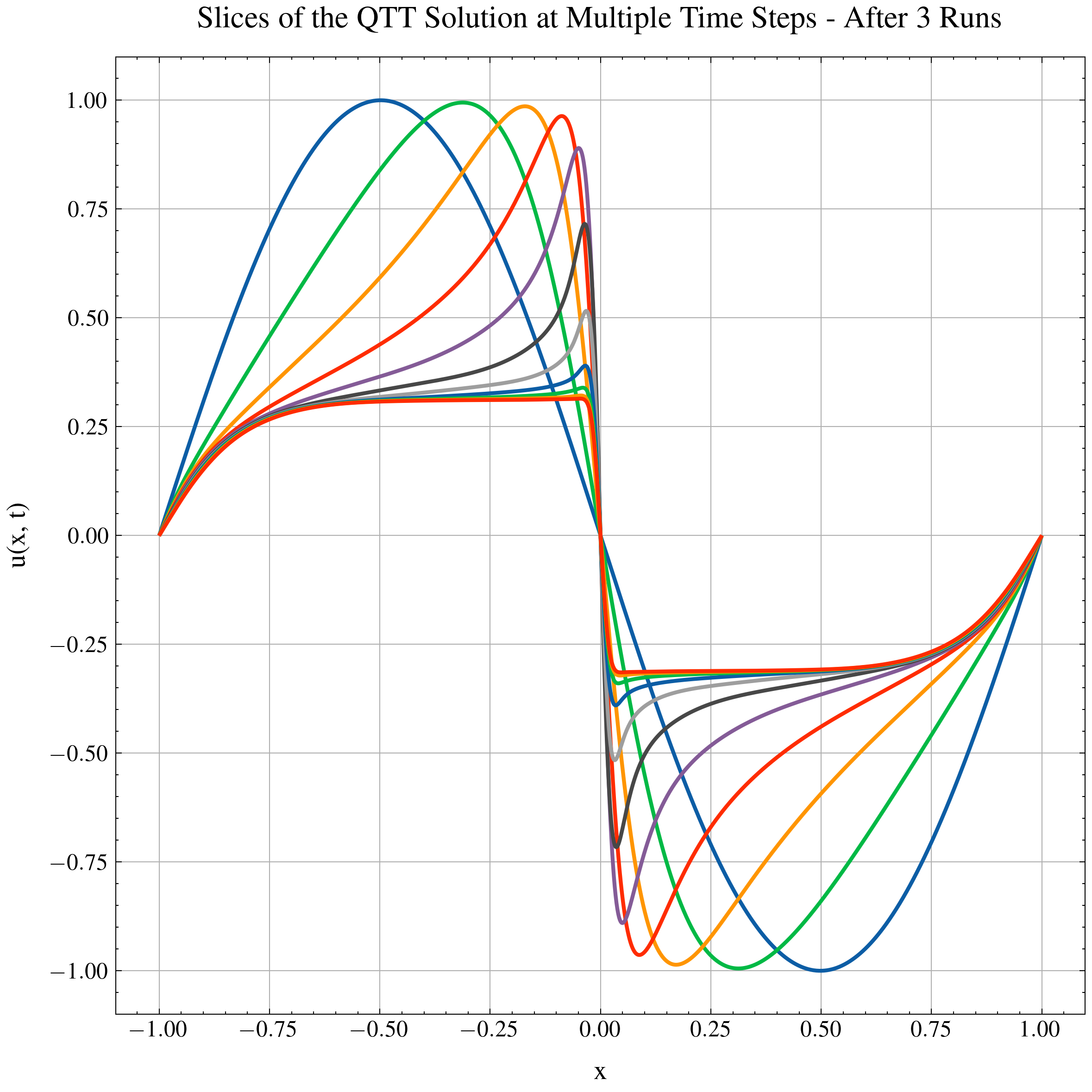}
    \end{minipage}
    \hfill
    \begin{minipage}{0.24\textwidth}
        \centering
        \includegraphics[width=\linewidth]{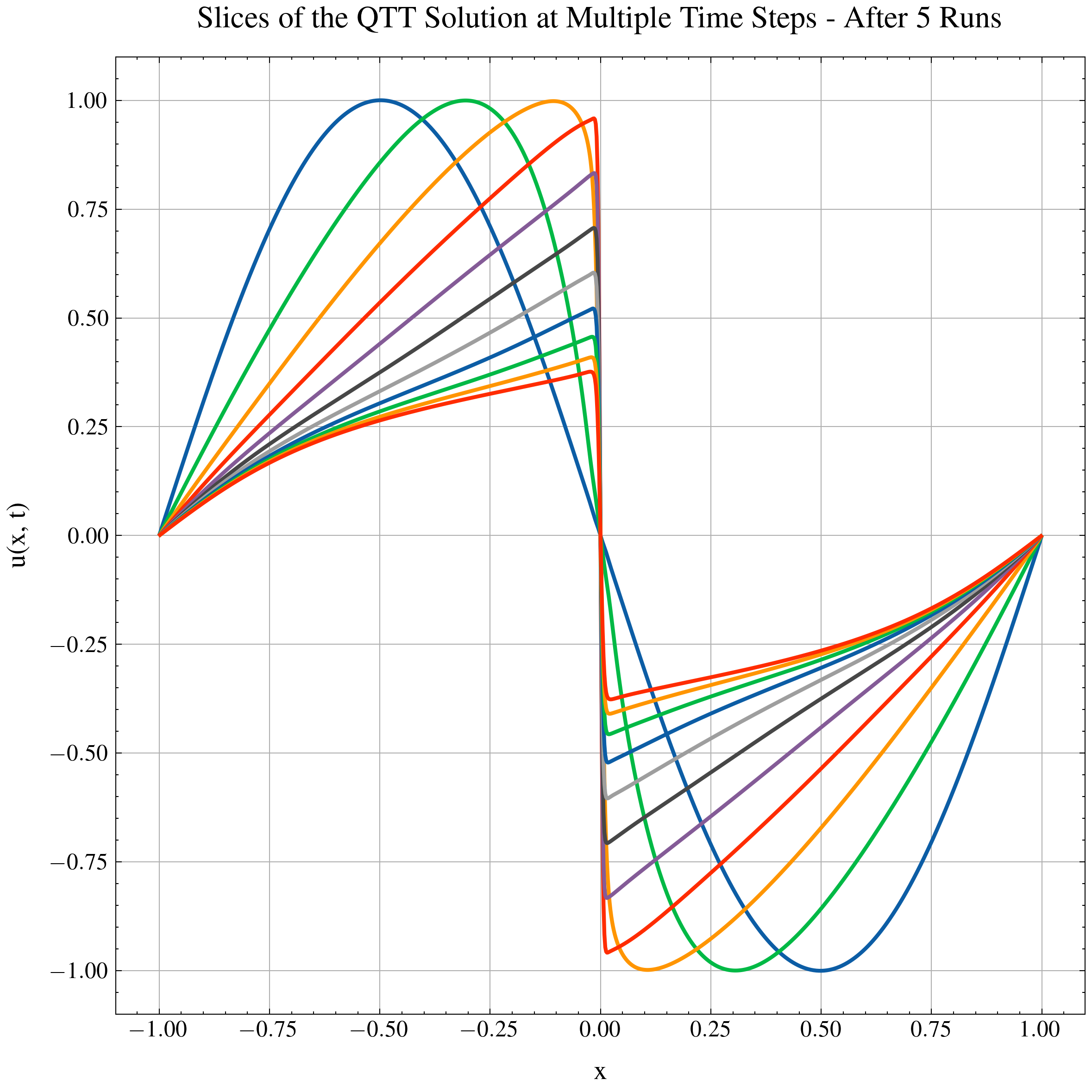}
    \end{minipage}
    \hfill
    \begin{minipage}{0.24\textwidth}
        \centering
        \includegraphics[width=\linewidth]{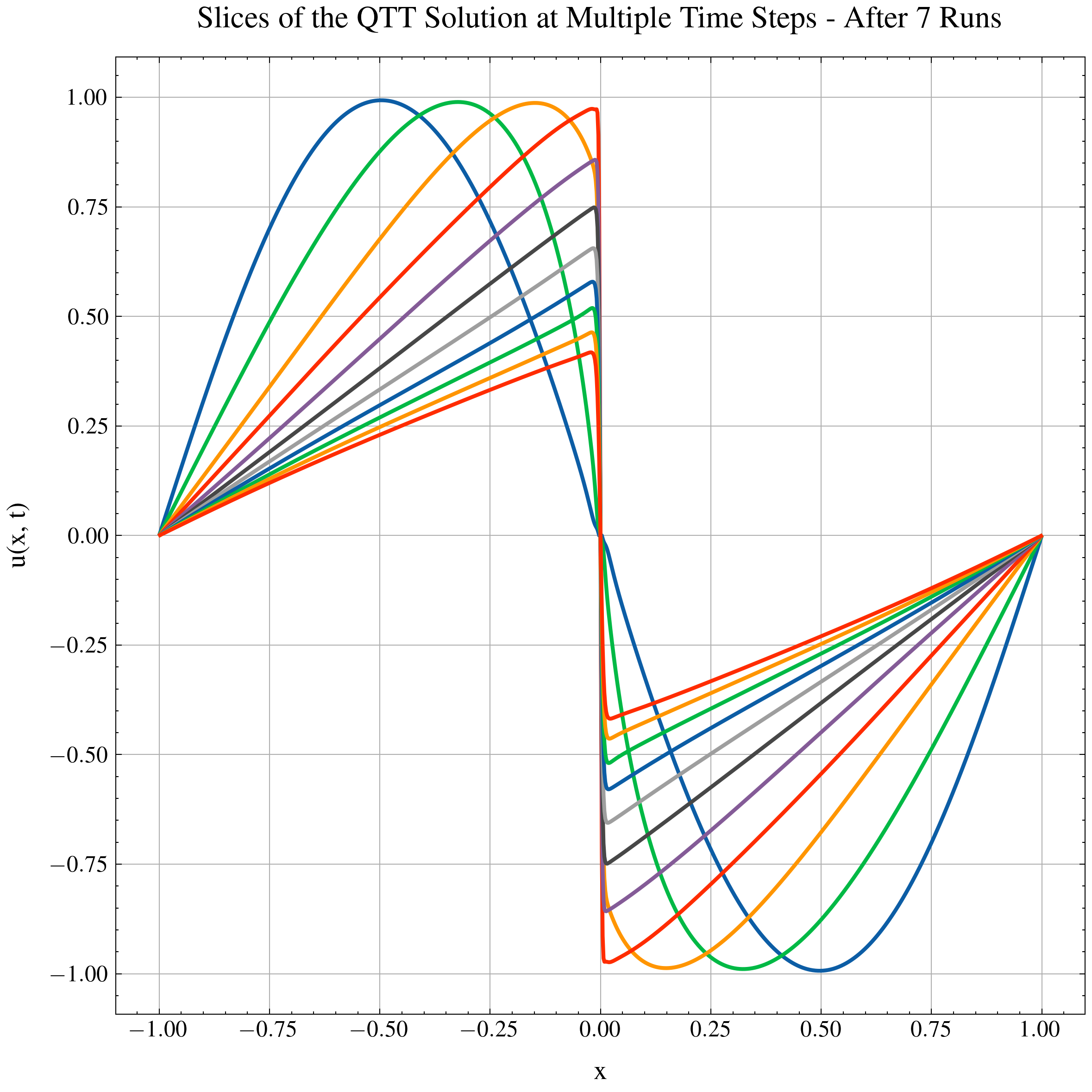}
    \end{minipage}

    \caption{Comparison of QTT solutions (top row) and slices (bottom row) for different numbers of \textit{runs} for 10 uniformed spaced time steps.}
    \label{fig:qtt_comparison_small}
\end{figure}
\end{center}
To demonstrate the convergence of our method, we present the leftmost table, which shows the resulting MSE between the solutions for different numbers of \textit{runs}. On the right, we compare the MSE of our 5 and 7-run QTT solution at different slices that correspond to six sequentially equally spaced time steps against an approximation of the analytical solution. This approximation, originally implemented in \cite{python_burgers_analytic}, is based on Hermite quadrature as described by Basdevant et al. in \cite{burgers_pinns_analytic}.

\begin{table}[h]
    \centering
    \begin{minipage}{0.45\textwidth}
        \centering
        \begin{tabular}{cc}
            \toprule
            \textbf{\textit{Runs} Compared} & \textbf{MSE} \\
            \hline
            2 vs. 3 & 0.094 \\
            3 vs. 5 & 0.0089 \\
            5 vs. 7 & 0.0005 \\
            7 vs. 8 & 1.34e-05 \\
            8 vs. 9 & 6.59e-07 \\
            \hline
        \end{tabular}
        \caption{MSE between the full solutions with different numbers of \textit{runs}.}
        \label{tab:mse_runs}
    \end{minipage}%
    \hfill
    \begin{minipage}{0.45\textwidth}
        \centering
        \begin{tabular}{ccc}
         \textbf{Equally Spaced} & \multicolumn{1}{c}{\textbf{MSE}} & \multicolumn{1}{c}{\textbf{MSE}} \\ 
        \cmidrule(lr){2-2}
        \textbf{Time Steps} & \textbf{5-run QTT}  & \textbf{7-run QTT} \\ 
        \hline
            $t_1$  & 6.81e-07 & 8.20e-09\\
            $t_2$  & 5.81e-04 & 1.34e-05\\
            $t_3$  & 5.58e-04 & 5.32e-06\\
            $t_4$  & 5.58e-04 & 5.31e-06\\
            $t_5$  & 5.81e-04 & 1.34e-05\\
            $t_6$  & 1.70e-07 & 2.05e-09\\
            \hline
        \end{tabular}
        \caption{MSE of the 5 and 7-run QTT solution against the analytical approximation based on Hermite quadrature at equally spaced time steps, starting from the first one after the initial condition.}
        \label{tab:mse_qtt_analytic}
    \end{minipage}
\end{table}
Our space-time QTT solver demonstrates rapid convergence to the solution, even when using a simple discretization scheme and a small number of of \textit{runs}. Despite its simplicity, our solver achieves high precision while maintaining an exceptionally low runtime—orders of magnitude faster than PINN-based approaches. The same instance of the Burgers' equation has been extensively benchmarked in \cite{pinns_bench} using various PINN configurations. According to their Tables 8 and 12, achieving the same average MSE as our 7-runs QTT solver requires approximately 284 seconds—nearly 100 times longer than our method for comparable accuracy.

~ \newpage
%\bibliography{main}
%\bibliographystyle{abbrvnat}

\end{document}